\documentclass[a4paper, 10pt]{amsart}

\usepackage{mathptmx}
\usepackage[latin1]{inputenc}

\usepackage[alphabetic]{amsrefs}
\usepackage{graphicx}
\usepackage{amstext}
\usepackage{amsfonts}
\usepackage{amssymb}
\usepackage{amscd}
\usepackage{amsmath}
\usepackage{amsthm}
\usepackage{setspace}
\usepackage{graphicx}
\usepackage[all]{xy}
\usepackage{wasysym}
\usepackage{multicol}

\usepackage{acronym}

\newcommand {\F}{\mathcal{F}}

\newcommand {\C}{\mathbb{C}}
\newcommand {\R}{\mathbb{R}}
\newcommand {\E}{\mathbb{E}}
\newcommand {\Rn}{(\mathbb{R}^n, \langle \cdot, \cdot \rangle)}

\renewcommand {\H}{\vec{H}}
\newcommand {\tnu}{\tilde{\nu}}
\newcommand {\tr}{\mathrm{tr}}

\theoremstyle{plain}
\newtheorem{thm}{Theorem}[section]
\newtheorem{lem}[thm]{Lemma}
\newtheorem{prop}[thm]{Proposition}
\newtheorem{cor}[thm]{Corollary}
\theoremstyle{definition}
\newtheorem{defn}[thm]{Definition}

\theoremstyle{remark}
\newtheorem{rem}[thm]{Remark}

\begin{document}

\title[Self-Shrinkers of the MCF] 
{Spacelike Self-Similar Shrinking Solutions of the Mean Curvature Flow in Pseudo-Euclidean Spaces}

%
\author{M\'arcio Rostirolla Adames}
\address[M\'arcio Rostirolla Adames]{Departamento Acad\^emico de Matem\'atica\\
        UTFPR - Campus Curitiba\\
        Curitiba - PR, Brasil 80230-901}
\email[M.~Adames]{marcioadames@utfpr.edu.br}
\thanks{This work was done with financial support of the CNPq - National Counsel of Scientific and Technologic Development from Brazil - and developed at the \textit{Leibniz Universit\:at Hannover} during the author's Ph. D.}

\date{July 8, 2013}


\maketitle

\begin{abstract}
I classify spacelike self-similar shrinking solutions of the mean curvature flow in pseudo-euclidean space in arbitrary codimension,
if the mean curvature vector is not a null vector and the principal normal vector is parallel in the normal bundle.
Moreover, I exclude the existence of such self-shrinkers in several cases. 
The classification is analogous to the existing classification in the euclidean case \cites{Hui93, SMK05}.
\end{abstract}

%
%
%
%
%

\section{Introduction}
The Mean Curvature Flow (MCF) of an immersion $F: M \to N$ of a smooth manifold $M$ into a Riemannian manifold $(N,h)$ is a natural way to deform this immersion into something ``rounder'' or ``more regular''. It is a smooth family of isometric immersions $F_t: M \to N$, $t \in [0,T)$ that satisfies
$$
\frac{d F_t}{dt} = \H, \hspace{15mm} F_0(x) = F(x).
$$

The MCF has been studied by many. It not ony regularizes the initial surface but also produces singularities. Suppose now that the target manifold $N$ is the Euclidean space $\E^n$. In an important work on the MCF of convex compact hypersurfaces \cite{Hui84} Huisken showed, among other results, that the supremum of the norm of the second fundamental form $\sup_M \|A\|^2$ explodes as $t \to T$ (the maximal existence time) if there is a finite time ($T < \infty$) singularity. This happens because an upper bound on the second fundamental form would imply upper bounds on all the derivatives $\nabla^{(k)} A_{ij}$ and the solution could be then extended beyond $T$, which 
is a contradiction. However, this can be done not only for hypersurfaces, but for a broader class of manifolds and in any codimension (see \cite{SMK11} Proposition 3.11 and Remark 3.12 or \cite{MR2801634}).

In a subsequent work \cite{Hui90}, Huisken showed, with his famous monotonicity formula, that hypersurfaces satisfying a natural\footnote{This is the growth rate of some simple hypersurfaces, like spheres and cylinders.} growth in the norm of the second fundamental form:
$$
\max_{M} \|A\|^2 \leq \frac{C_0}{2(T-t)},
$$
for some constant $C_0 >0$, deform asymptotically near a singularity to self-similar solutions of the MCF after some blow up process (rescaling the surface and changing the time variable). This result depends only on the existence of some integrals with respect to the backwards heat kernel and holds, for example, if $M$ is closed. Later Ilmanen \cite{Il97} and White \cite{Wh94} proved that the all finite time singularities in the generalized sense of the Brakke flow \cite{Br} are self-similar solutions of the MCF.

These self-similar solutions of the MCF are also called self-shrinkers to avoid confusion with other types of solutions that preserve the ``form'' of the surface, like self-expanders and translating solutions. They are homotheties that shrink the initial manifold and are given by the equation
$$
\H = -F^\perp.
$$

Because of the relation between singularities of the MCF and self-shrinkers, there is interest in classifying and giving examples of these in special cases. Abresch and Langer \cite{A.L.86} gave the complete classification of the closed plane curves that shrink homothetically, they are the circles and the so called Abresch \& Langer curves. Huisken proved in \cite{Hui93} that the self-shrinking hypersurfaces with non-negative mean curvature (compact or non-compact) are spheres, cylinders and the product of an Abresch \& Langer curve with an affine space. The result of Huisken was later generalized by Smoczyk \cite{SMK05} for higher codimensional immersions, with the assumption that the principal normal is parallel in the normal bundle and $\|\H\|_{\E} \neq 0$. A related result was found by Cao and Li \cite{CaoLi} in any codimension: the self-shrinkers with $\|A\|^2 \leq 1$ are spheres, planes or cylinders. There are also Bernstein type results for self-shrinkers in higher codimension of Q. Ding and Z. Wang \cite{DingWang}, who generalize works of Lu Wang \cite{LuWang} and Huisken \& Ecker \cite{EcHu89}. Recently, Baker \cite{Baker} proved that high codimensional self-shrinkers under certain conditions for the second fundamental tensor and mean curvature vector are spheres or cylinders.

For hypersurfaces in $\R^3$, there are examples of a shrinking doughnut of Angenent \cite{Ang} and many numerical examples of Chopp \cite{Ch} and Ilmanen \cite{Il98}, like ``punctured saddles`` made of many handles crushing at the same time, which are highly unstable, depending on the surface having many symmetries. Colding and Minicozzi \cite{ColMin} showed that the only stable singularities for smooth closed embedded surfaces in $\R^3$ are cylinders and spheres. For the Lagrangian MCF, Joyce, Lee and Tsui \cite{JLT}, Anciaux \cite{Anc} and Wang \cite{Wa} have examples. There are other results in different contexts.

The main purpose of this work is to study self-shrinkers of the MCF in higher codimension in the pseudo-euclidean case. By that we mean that the target manifold $N$ is a pseudo-euclidean space, so that the most interesting new case is the Minkowski space $\R^{1,n}$. The MCF of spacelike hypersurfaces in the Minkowski space was studied for example by Ecker \cite{Ecker97} and a related flow was considered by Ecker and Huisken \cite{EcHu}. Gerhardt \cite{MR2488947} also studies curvature flows in semi-Riemannian manifolds, specially the inverse mean curvature flow. Beyond this Bergner and Sch\"afer \cite{BerSchae} considered the mean curvature flow in the 3-dimensional Minkowski space. Furthermore Li \& Salavessa have some results for the MCF of spacelike graphs \cite{MR2523502} in product manifolds.

In the first section are given some fundamental equations in order to fix the notation and in the second one we consider homotheties of the MCF that lie in hyperquadrics, and find,
similarly to Smoczyk's result in \cite{SMK05} for spheres in the Euclidean space $\E^n$, that the homotheties (with nondegenerate first fundamental form) of the MCF with initial immersion contained in a hyperquadric are exactly the minimally immersed submanifolds of the hyperquadric if $k>0$ or $k<0$, as Theorems \ref{PositiveSphere} and \ref{NegativeSphere} state. Moreover, given the initial minimal immersion, the flow can be explicitly calculated. If $k = 0$ (the light cone), a homothety with nondegenerate first fundamental form would immediately leave the light cone and thence could not be a homothety starting at $t=0$ because the light cone is star shaped, as stated in Theorem \ref{nullsphere}. But, as Ecker noted in \cite{Ecker97}, the upper light-cone would immediately change to a 
hyperquadric and the explicit solution of the MCF with the upper light-cone as initial condition in the Minkowski space ($\R^{1,n}$) would be the graph of the function $\delta (x,t) = \sqrt{\|x\|_{\E}^2 + 2(n-1)t}$, which is a homothety after $t=0$.

There is a big difference between the flow of minimal surfaces of the hyperquadrics with $k>0$ and the ones with $k<0$; if $k>0$, they shrink to a point (at least the compact ones) in finite time, but if $k < 0$, they expand and never produce singularities. Beyond this, they are given by different equations. The following results are for the shrinking\footnote{The expanding case satisfies $\H = F^\perp$.} case, which are the isometric immersions $F:M \to \Rn$ satisfying
$$
\H=- F^\perp.
$$

Our ``domain'' manifold $M$ is always assumed to be smooth, path connected, complete and orientable.

If one considers the self shrinkers and self-expanders that are contained in the hyperquadrics as submanifolds in the pseudo-euclidean space $(\R^n, \langle \cdot, \cdot \rangle)$, then one observes that $\nabla^\perp \H \equiv 0$ and $\nabla^\perp \nu \equiv 0$, where $\nu:= \H/\|\|\H\|$ is the principal normal. A natural question is whether these conditions are also sufficient to guarantee that a spacelike\footnote{We use elliptic methods to obtain our results (maximum principles, that do not hold for hyperbolic equations).} self-shrinker lies in a hyperquadric. The condition $\nabla^\perp \H = 0$ implies this immediately if $M$ is compact, because $\|\H\|^2$ is then constant and the maximum principle implies, with equation
\begin{equation}\label{negociodoidoloco}
\triangle \|F\|^2 = 2m - 2\|\H\|^2,
\end{equation}
that $\|F\|^2$ is constant. So, in this work, we examine the self-shrinkers of the MCF with $\|\H\|^2 \neq 0$ and $\nabla^\perp \nu = 0$. The condition $\nabla^\perp \nu = 0$ is natural because it holds for any hypersurface.

The third section deals with fundamental equations for self-shrinkers with the principal normal parallel in the normal bundle and the compact case. The fourth section is about the non-compact case.

Equation \eqref{negociodoidoloco} already shows that there are no compact self-shrinkers with $\|\H\|^2<0$. In this article, the inexistence of self-shrinkers with $\|\H\|^2<0$ is proven, also in the non-compact case under certain hypothesis, as stated in Theorem \ref{SumUpNoSS}.



As a consequence of this, the Minkowski space does not (in all of our treated cases) have spacelike self-shrinking hypersurfaces. This could already be seen from Ecker's longtime existence result for spacelike hypersurfaces in Minkowsky space \cite{Ecker97}.

If $\|\H\|^2 > 0$ one finds, just as Smoczyk in \cite{SMK05} for the Euclidean case, that if $M$ is compact and dim($M$) $\geq 2$, the only spacelike self-shrinkers of the MCF with $\|\H\|^2(p) \neq 0$, $\forall p \in M$, and $\nabla^\perp \nu \equiv 0$ are the minimal\footnote{By minimal we mean the ones satisfying $\H=0$. We use this name because the condition $\H=0$ is then mnemonic, although this condition does not imply minimality of the volume functional in pseudo-euclidean spaces, so they are just critical points of the volume functional.} submanifolds of hyperquadrics, as stated in Theorem \ref{PrinParNorSphe}.
%
This dimensional restriction is in fact optimal because in dimension one there are the Abresch \& Langer curves which are self-shrinkers and do not lie in spheres.


Again following Smoczyk in the non-compact case, one finds that the self-shrinkers with $\|\H\|^2 >0$ are products of affine spaces with minimal submanifolds of hyperquadrics or with homothetic solutions of the curve shortening flow\footnote{The curve shortening flow is the MCF for plane curves.} as stated in Theorem \ref{ResultNoncompact3}.


The proof of these Theorems is long and internally divided in lemmas to make its several steps easier to recognize. It was necessary to divide the proof in two cases. In both of them we split $TM$ into two involutive distributions. Then we use the Theorem of Frobenius to get foliations on $M$ whose leaves are totally geodesic immersed in $M$. After this, we calculate a formula that relates the second fundamental tensor of $F$ with these distributions. In particular the second fundamental tensor of $F$ is zero when restricted to one of these distributions, so that the leaves of this distribution are totally geodesic in $(\R^{q,n}, \langle \cdot, \cdot \rangle)$ and then, considering parallel transports inside these leaves, one finds that they are parallel affine subspaces of $\R^{q,n}$. The other distribution delivers the $\mathcal{H}_r$ and $\Gamma$ parts in the last Theorem. We get this considering the second fundamental tensor and mean curvature vectors of the inclusion of the leaves related 
to this distribution, with some extra effort to prove that $\Gamma$ lies on a plane (based on an idea of \cite{Hui93}). In the last step we construct an explicit map from these second leaves times $\R^{m-r}$ onto $F(M)$.

\begin{section}{Geometric Background}
Let $(\R^n, \langle \cdot, \cdot \rangle)$ be an inner product space. This means that $\langle \cdot, \cdot \rangle$ is a symmetric bilinear form which is nondegenerate (but not necessarily positive definite).

\begin{defn}
For $n \in \{2,3, \ldots\}$ we call the set
\begin{displaymath}
\mathcal{H}^{n-1}(k):= \{x \in \R^n : \|x\|^2 = k \}
\end{displaymath}
\textit{Hyperquadric} of dimension $n-1$ and parameter $k$, $k \in \R$ fixed.
\end{defn}

%
%
%
%
%
%
%
%
%
%
%
%
%

 Let $F: M \to \Rn$ an immersion. The identification of $T_p\R^n$ with $\R^n$ induces a semi-Riemannian metric (denoted also by $\langle \cdot, \cdot \rangle$) on $T_p\R^n$ and the immersion F induces a semi-Riemannian metric $g:= F^* \langle \cdot, \cdot \rangle$ over $M$, if it is nondegenerate. We assume that $g$ is nondegenerate. Let $\nabla^g$ be the Levi-Civita connection induced by $g$. 
 Then:
$$
dF(\nabla^g_X Y) = (D_{dF(X)} dF(Y))^\top.
$$

We also use the following connections on several bundles:
\begin{itemize}
\item $\nabla^{F^*T\R^n}$ on the pullback bundle defined as $\nabla^{F^*T\R^n}_X Y := D_{dF(X)} Y$ for any $X \in \Gamma(TM)$ and $Y \in \Gamma(F^*T\R^n)$.
\item $\nabla^\perp$ on the normal bundle defined as $\nabla^\perp_X Y = (D_{dF(X)} Y)^\perp$ for any $X \in \Gamma(TM)$ and $Y \in \Gamma(TM^\perp)$.
\item $\nabla^*$ on the dual of a bundle $E$ over $M$ defined through $(\nabla^*_X \epsilon)(e) :=X(\epsilon(e)) - \epsilon (\nabla_X e)$ for any $X \in \Gamma(TM)$, $e \in \Gamma(E)$ and $\epsilon \in \Gamma(E^*)$.
\item $\nabla^{E \otimes F}$ on the product bundle $E \otimes F$ of two bundles $E$ and $F$ over $M$ defined as $\nabla^{E \otimes F}_X (e\otimes f):= \nabla^E_X e \otimes f + e\otimes \nabla^F_X f$
\end{itemize}
We usually omit most of the superscript indicating the bundle. We just use $\nabla$ for most of the cases and $\nabla^\perp$ if we project, on the normal bundle, the component of the tensor that lies in $TM^\perp$. For example, for $X \in \Gamma(TM)$ and $Y \otimes Z \in \Gamma(TM \otimes TM^\perp)$, it holds
\begin{align*}
\nabla_X (Y\otimes Z) =& (\nabla_X Y)\otimes Z + Y \otimes (\nabla_X Z)\\
\nabla_X^\perp (Y\otimes Z) =& (\nabla_X Y)\otimes Z + Y \otimes (\nabla_X^\perp Z).
\end{align*}
%
%
%
%
%

\begin{rem}
We write $A$ and $\H$ (sometimes $A_F$ or $\H_F$) for the second fundamental tensor and the mean curvature vector of an isometric immersion $F$. 
\end{rem}

%

We use Latin letters for indices of tensors on $M$ and Greek letters for indices of tensors on the target manifold N, in our case $\Rn$. We also use the Einstein's convention for sums. So that:
\begin{align}
\nonumber A_{ij}&:= (\nabla dF)\left(\frac{\partial}{\partial x^i}, \frac{\partial}{\partial x^j}\right) = \left(\nabla \left(\frac{\partial F^\alpha}{\partial x^k} dx^k \otimes \frac{\partial}{\partial y^\alpha} \right)\right) \left(\frac{\partial}{\partial x^i}, \frac{\partial}{\partial x^j}\right) \\
\nonumber &= \frac{\partial^2 F^\alpha}{\partial x^i \partial x^j} \frac{\partial}{\partial y^\alpha} + \frac{\partial F^\alpha}{\partial x^k} (\nabla dx^k) \left(\frac{\partial}{\partial x^i}, \frac{\partial}{\partial x^j}\right)\otimes \frac{\partial}{\partial y^\alpha}\\
\nonumber &= \frac{\partial^2 F^\alpha}{\partial x^i \partial x^j} \frac{\partial}{\partial y^\alpha} - \frac{\partial F^\alpha}{\partial x^k} \left(dx^k\left(\nabla_{\frac{\partial}{\partial x^i}} \frac{\partial}{\partial x^j} \right) \right)\otimes \frac{\partial}{\partial y^\alpha}\\
\label{basicSecondFunForm} &= \frac{\partial^2 F^\alpha}{\partial x^i \partial x^j} \frac{\partial}{\partial y^\alpha} - \frac{\partial F^\alpha}{\partial x^k} \Gamma^k_{ij} \frac{\partial}{\partial y^\alpha} = \nabla_i \nabla_j F,
\end{align}
considering $F \in \Gamma(F^*T\R^n)$.

We use the (rough) Laplacian $\triangle$ on sections of several bundles,  and write $\triangle := g^{ij} \nabla_i \nabla_j$ and $\triangle^\perp := g^{ij} \nabla_i^\perp \nabla_j^\perp$. Beyond this the second fundamental tensor is written $A_{ij}= \nabla_i \nabla_j F$ and it follows $\H = \triangle F$.

We use following conventions for the Riemannian curvature vector of the tangent bundle in local coordinates:
$$
R^l_{\,\,kij} \frac{\partial}{\partial x^l}= R\left(\frac{\partial}{\partial x^i},\frac{\partial}{\partial x^j}\right)\frac{\partial}{\partial x^k},
$$
and
$$
R_{skij}=R\left(\frac{\partial}{\partial x^s}, \frac{\partial}{\partial x^k}, \frac{\partial}{\partial x^i}, \frac{\partial}{\partial x^j}\right) = R^l_{\,\,kij} g_{ls}.
$$

The Codazzi equation in local coordinates as
\begin{equation}
 \label{Codazzi Equation} \nabla_l A_{ij} - \nabla_i A_{lj} = R^k_{\,\,jli}F^\alpha_k
\end{equation}
and considering $A$ as a section in the normal bundle, $A \in \Gamma(TM^\perp \otimes TM^* \otimes TM^*)$,
\begin{align}
\nonumber \nabla_l^\perp A_{ij} - \nabla_i^\perp A_{lj} &= [\nabla_l A_{ij} - \nabla_i A_{lj}]^\perp\\
\label{codazzi}\nabla_l^\perp A_{ij} - \nabla_i^\perp A_{lj} &=[R^k_{jli}F^\alpha_k]^\perp =0.
\end{align}


We make use of Gau{{\ss}} equation:
\begin{equation}\label{gauss}
R_{klij} = \langle A_{ik}, A_{jl} \rangle - \langle A_{jk}, A_{il} \rangle.
\end{equation}
and the Ricci-equation:
\begin{equation}
R^\perp (X,Y) \eta = \tr (\left< \eta, A\left(Y,\cdot \right) \right> A\left(X, \cdot \right)) -\tr(\left< \eta, A\left(X,\cdot\right) \right> A\left(Y,\cdot\right)).
\end{equation}

The Riemannian curvature tensor of the normal bundle $R^\perp_{ij}$ can be seen as the section $\langle R^\perp \left(\frac{\partial}{\partial x^i},\frac{\partial}{\partial x^j} \right) \cdot, \cdot \rangle \in \Gamma(T^*M^\perp \otimes T^*M^\perp)$. 
The Ricci equation is then written
\begin{equation}\label{RicciWedge}
R^\perp_{ij} = A_{jk}\otimes A^k_i - A_{ik}\otimes A^k_j =: A_{jk}\wedge A^k_i.
\end{equation}

Furthermore we need the commutation formula:
\begin{lem}\label{commutation}
Let $M$ be a $m$-dimensional differentiable manifold, $E$ be a vector bundle of dimension $\delta$ over $M$ and $\nabla$ be a connection on this bundle. If $T^{\alpha_1 \ldots \alpha_\phi}_{k_1 \ldots k_r} \in \Gamma(T^*M \otimes \ldots \otimes T^*M \otimes E \otimes \ldots \otimes E)$ then
\begin{align}
\nabla_i \nabla_j T^{\alpha_1 \ldots \alpha_\phi}_{k_1 \ldots k_r} - \nabla_j \nabla_i T^{\alpha_1 \ldots \alpha_\phi}_{k_1 \ldots k_r} =& - \sum_{h=1}^r \sum_{p=1}^m R^p_{k_hij}T^{\alpha_1 \ldots \alpha_\phi}_{k_1 \ldots k_{h-1}p k_{h+1} \ldots k_r}\\
\nonumber &+ \sum_{h=1}^\phi \sum_{\beta=1}^\delta R^{\alpha_h}_{\beta ij}T_{k_1 \ldots k_r}^{\alpha_1 \ldots \alpha_{h-1}\beta \alpha_{h+1} \ldots \alpha_\phi}.
\end{align}
\end{lem}

\end{section}

\begin{section}{Hyperquadric Homotheties of the MCF}

\begin{subsection}{Hyperquadric Homotheties of the MCF}


\begin{defn}
Let $M$ be a smooth manifold, $(N,h)$ be a semi-Riemannian manifold and $F_0:M \to N$ be an immersion. A smooth a family of isometric immersions $F:M \times [0,T) \to N$, for some $T>0$, such that the metric $g_t:=F(\cdot, t)^* h$ is nondegenerate for all $t \in [0,T)$ is called a \textit{solution of the mean curvature flow} with initial immersion $F_0$ if it satisfies
\begin{equation}
\frac{dF}{dt}(p,t)=\H(p), \hspace{10mm} \mbox{and} \hspace{10mm} F(p, 0) = F_0(p) \,\, \forall p \in M, t \in [0,T),
\end{equation}
where $\H$ is the mean curvature vector of the immersion $F(\cdot, t): M \to (N,h)$.
\end{defn}

Now we consider properties of homotheties generated by the mean curvature flow. Let $F: M \times [0,T) \to \Rn$ be a solution of the MCF for some initial immersion, such that there is a rescaling function $c:[0,T) \to (0, \infty)$, with $c(0)=1$, so that $\tilde{F}:=cF$ satisfies
\begin{equation}\label{originselfsimilar}
\left< \frac{d\tilde{F}}{dt}(p,t),\tnu \right> = 0 \hspace{25mm} \forall \, \tnu \in T_pM ^{\perp},
\end{equation}
which implies
\begin{equation}\label{selfie}
\H=-\frac{\dot{c}}{c}F^\perp.
\end{equation}

\begin{defn}
Let $F: M \times [0,T) \to \Rn$ be a solution of the MCF for some initial immersion. If there is a differentiable function $c:[0,T) \to (0, \infty)$ with $c(0)=1$ such that
\begin{equation}\label{selfsimilar}
\H=-\frac{\dot{c}}{c}F^\perp,
\end{equation}
we say that $F$ is a homothety of the MCF.
\end{defn}



In particular we look now at the hyperquadric homotheties of the MCF in $\Rn$


From now on let $F:M \times [0,T) \to \Rn$ be a homothety of the MCF. It could happen that some solutions of the flow in which the initial immersion lies in a hyperquadric, i. e. $\|F(x,0)\|=k$ for all $x \in M$, cease lying in some hyperquadric during the flow. This cannot happen for homotheties, as the following result states:

\begin{lem}\label{primeirolema}
If $F(0,x) \subset \mathcal{H}^{n-1}(k(0))$ for all $x \in M$ then $F(t,x) \subset \mathcal{H}^{n-1}(k(t))$ for all $x \in M$, for some function $k:[0,T) \to \R$.
\end{lem}


As the position vector in a hyperquadric is normal (with respect to the inner product that generates the hyperquadric) to the hyperquadric, it follows that $F$ is always orthogonal to $\mathcal{H}^{n-1}(k(t))$.With this, $\|F\|^2$ can be calculated:


\begin{lem}\label{keepspherical}
$\|F(t)\|^2 = k(0) - 2mt$ for $t \in [0,T)$.
\end{lem}
\begin{proof}[Proof]
For $t=0$ it is clear that $\|F(0)\|^2 = k(0)$. And for all $t \in [0,T)$:
\begin{align*}
\frac{d}{dt} \|F(t)\|^2 &= 2 \left< \frac{d}{dt} F(t), F(t)\right> = 2\left< \vec{H}(t), F(t)\right>\\
&= 2\left< \nabla^j\nabla_j F(t), F(t)\right>= -2\left< F(t)_j, F(t)_l \right> g^{jl} = -2m.
\end{align*}
\end{proof}

%

We prove now that a hyperquadric homothety of the mean curvature flow is a minimal immersion in the hyperquadric $\mathcal{H}^{n-1}(k(0)-2mt)$ for all $t \in [0,T)$. We will need the following Lemma:

\begin{lem}\label{reletionsbetweenHs}
Let $F:M \to N$ and $G:N \to P$ be isometric immersions between semi-Riemannian manifolds $(M,g),(N,h)$ and $(P,l)$. Denote $\H_F$, $\H_G$ and $\H_{G\circ F}$ the mean curvatures of $F$, $G$ and $G\circ F$ respectively. Then:
\begin{equation}
 \label{chainrulesecondfund}(\nabla d(G\circ F))_x = (\nabla dG)_{F(x)} (dF\cdot, dF\cdot) + dG_{F(x)}\circ \nabla dF
\end{equation}
and $\H_{G\circ F} = dG(\H_F) + \tr_M(\nabla dG)\left(dF\cdot, dF\cdot \right)$.
\end{lem}
%
%
\begin{thm}\label{self-similar is minimal}
 Let $F: M \times [0,T) \to \Rn$ be a hyperquadric homothety of the mean curvature flow, then $F(M,t)$ is a minimal surface of the hyperquadric $\mathcal{H}^{n-1}(\|F(0)\|^2 -2m t)$ for all $t \in [0,T)$.
\end{thm}
\begin{proof}
We consider the natural inclusion $I(t)$ of the hyperquadric $\mathcal{H}^{n-1}(k(0) - 2mt)$ into $\Rn$ 
and the immersion $G:M \to \mathcal{H}^{n-1}(k)$ defined as $G:=I^{-1} \circ F$, as in the diagram:
\begin{displaymath}
    \xymatrix{\mathcal{H}^{n-1}(k(0) - 2mt) \ar[r]^{\,\,\,\,\,\,\,\,\,\,\,\,\,I} & \Rn \\
        M \ar[u]^G \ar[ur]_F & }
\end{displaymath}
Writing $\H_F$, $\H_G$ $\H_I$ for the respective mean curvature vectors, it holds:
\begin{itemize}
 \item $\H_F \in TM^\perp$,
 \item $dI(\H_G) \in (dI(T\mathcal{H}^{n-1}))$,
 \item $g^{ij}\left(\nabla dI\right)\left(dG\left(\frac{\partial}{\partial x^i}\right), dG\left(\frac{\partial}{\partial x^j}\right) \right) \in \left.T\mathcal{H}^{n-1}\right.^\perp$.
\end{itemize}
But equation \eqref{selfsimilar} implies that $\H_F \in (T\mathcal{H}^{n-1})^\perp$. Thus, $dI(\H_G)$ is the only term tangential to the hyperquadric in Lemma \ref{reletionsbetweenHs}, thence $dI(\H_G)=0$ and $\H_G =0$ (for $I$ is an immersion).
\end{proof}

Further, we can calculate $-\frac{\dot{c}}{c}$. Let $t \in [0,T)$ be fixed and $x \in M$ be any point
\begin{equation}
\nonumber -\frac{\dot{c}}{c} \|F(x,t)\|^2 = \langle \H_F, F \rangle = -\langle F_i, F_j \rangle g^{ij} = -m \Rightarrow -\frac{\dot{c}}{c}=-\frac{m}{\|F(t)\|^2},
\end{equation}
and Lemma \ref{keepspherical} implies that
\begin{equation}\label{meancurvaturesphere}
 \H_F(t) = -\frac{m}{\|F(0)\|^2 -2mt}F(t).
\end{equation}
\end{subsection}

\begin{subsection}{Existence and Uniqueness}
\begin{subsubsection}{Immersion in the Hyperquadric $\mathcal{H}^{n-1}(k)$ with $k>0$.}

Let $\|F(x,0)\|^2 = k > 0$ for all $x \in M$. The rescaling function $c(t)$ is given, from eq. \eqref{selfsimilar} and \eqref{meancurvaturesphere}, by $c(t):= \sqrt{k}(k -2mt)^{-1/2}$.

It follows from equation \eqref{meancurvaturesphere}, for any $(x,t) \in M \times [0,T)$, that
\begin{equation}\label{HSphere}
\frac{d}{dt}F(x,t) = \vec{H}_{F(\cdot,t)}(x) = - \frac{m}{k - 2mt}F(x,t) = -\frac{\dot{c}}{c}F(x,t) \Longrightarrow \frac{d}{dt}(cF(t)) = 0.
\end{equation}

Hence $cF(x,t) = F(x,0)$ and
\begin{equation}\label{explicit}
 F(x,t) = \frac{1}{c}F(x,0).
\end{equation}
By construction we proved that if there is a hyperquadric homothety of the MCF, then it has to be given by eq. \eqref{explicit}. So the solution is unique in the class of hyperquadric homothetic solutions.
We still have to deal with the question of existence. As in Theorem \ref{self-similar is minimal}, a hyperquadric homothety of the mean curvature flow has to be a minimal surface of the hyperquadric. This motivates the following Theorem:

\begin{thm}\label{PositiveSphere}
Let $F: M^m \to \Rn$ be an immersion such that $g:=F^*\langle \cdot, \cdot \rangle$ is nondegenerate and $\|F \|^2 = k \in \R$, $k>0$, then the solution of the MCF of this initial immersion is a homothety if, and only if, $F:M \to \mathcal{H}^{n-1}(k)$ is a minimal immersion in the hyperquadric $\mathcal{H}^{n-1}(k)$. The mean curvature flow of $F$ has a solution $F: M \times [0 ,T) \to \Rn$ with $T=\frac{k}{2m}$; moreover, the solution is $F(x,t):=c^{-1}(t) F(x)$, with $c(t):= \sqrt{k}(k -2mt)^{-1/2}$, $\forall (x,t) \in M \times [0 ,T)$.
\end{thm}
\begin{proof}
We have to prove that the homothety given by eq. \eqref{explicit} is a solution of the mean curvature flow. Let us write $F(t):=F(\cdot,t)$ and $I$ for the inclusion of $\mathcal{H}^{n-1}(k)$ into $\Rn$ and $G:=I^{-1} \circ F$. By Lemma \ref{reletionsbetweenHs} it follows
$$
\H_{F(0)} =g^{ij}(\nabla dI)\left(dG \left(\frac{\partial}{\partial x^i} \right), dG \left(\frac{\partial}{\partial x^j} \right)\right),
$$
because $F(0)$ is a minimal immersion on the hyperquadric. Moreover, $\H_{F(0)}$ is orthogonal to $\mathcal{H}^{n-1}(\|F(0)\|^2)$, but so is $F(0)$, such that there is a function $\varphi: M \to \R$ with $\H_{F(0)}=\varphi F(0)$. One can calculate $\varphi$:
\begin{align*}
\varphi \|F(0)\|^2 &= \langle \H_{F(0)} , F(0) \rangle = -g^{ij} \langle \nabla_j F(0), \nabla_i F(0) \rangle = -m\\
&\Longrightarrow \H_{F(0)} = -\frac{m}{\|F(0)\|^2}F(0).
\end{align*}
Now consider $F(t)=c^{-1}(t)F(0)$ (as in eq. \eqref{explicit}). Then $g^{ij}(t)=c^2(t)g^{ij}(0)$ and
$$
\H_{F(t)} = c(t) \H_{F(0)} = -\frac{m}{\|F(0)\|^2} c(t) F(0).
$$
On the other hand, for the function $c=\sqrt{k}(k -2mt)^{-1/2}$,
$$
\frac{dF(t)}{dt} = \frac{d}{dt}\left(\frac{1}{c(t)}\right)F(0)= -\frac{m}{\|F(0)\|^2}c(t) F(0) = \H_{F(t)}.
$$
Therefore this is a solution of the mean curvature flow.
\end{proof}

\end{subsubsection}

\begin{subsubsection}{Immersion in the Hyperquadric $\mathcal{H}^{n-1}(k)$ with $k<0$.}

Let $\|F(x,0)\|^2 = k < 0$ for all $x \in M$. The rescaling function $c(t)$ is given, from eq. \eqref{selfsimilar} and \eqref{meancurvaturesphere}, by $c(t):= \sqrt{-k}(-k+2mt)^{-1/2}$.

It follows from equation \eqref{meancurvaturesphere}, for any $(x,t) \in M \times [0,T)$, that
\begin{equation}\label{HSphere1}
\frac{d}{dt}F(x, t) = \vec{H}_F(\cdot, t)(x) = - \frac{m}{k - 2mt}F(x,t) = -\frac{\dot{c}}{c}F(x, t) \Longrightarrow \frac{d}{dt}(cF(x,t)) = 0.
\end{equation}

Hence is $cF(x,t) = F(x,0)$ and
\begin{equation}\label{explicit1}
 F(x,t) = \frac{1}{c}F(x,0).
\end{equation}
By construction we proved that if there is a hyperquadric homothety of the MCF, then it has to be given by eq. \eqref{explicit1}. So the solution is unique in the class of hyperquadric homothetic solutions. We still have to deal with the question of existence. This motivates the following Theorem:

\begin{thm}\label{NegativeSphere}
Let $F: M^m \to \Rn$ be an immersion such that $g:=F^*\langle \cdot, \cdot \rangle$ is nondegenerate and $\|F\|^2=k \in \R$, $k<0$, then the solution of the MCF of this initial immersion is a homothety if, and only if, $F:M \to \mathcal{H}^{n-1}(k)$ is a minimal immersion in the hyperquadric $\mathcal{H}^{n-1}(k)$. The mean curvature flow of $F$ has a solution $F(t): M \times [0 ,\infty) \to \Rn$; moreover, the solution is $F(x,t):=c^{-1}(t) F(x)$, with $c(t):=\sqrt{-k}(-k +2mt)^{-1/2}$, for all $(x,t) \in M \times [0,\infty)$.
\end{thm}
\begin{proof}
Analogous to Theorem \ref{PositiveSphere}.
\end{proof}

\end{subsubsection}

\begin{subsubsection}{Immersion in the Hyperquadric $\mathcal{H}^{n-1}(0)$}

Let $F:M \times [0,T)$ be a homothety generated by the MCF with $\|F(x,0)\|^2=0$ for all $x \in M$. From Lemma \ref{keepspherical} it holds $\|F(x,t)\|^2 = -2mt$ if $F^*\langle \cdot, \cdot \rangle$ is nondegenerate, so that 
\begin{equation}
\label{ijijijijijij}\|F(x,t)\|^2<0
\end{equation}
for all $(x,t) \in M \times (0,T)$.

On the other hand, $c(t)F(x,t) \in F(M,0)$ because $F$ is a homothety, so that
$$
0 = \|c(t)F(x,t)\|^2 = c(t)^2\|F(x,t)\|^2
$$
But $c(t)\neq 0$ because $F(M,t)=\{0\}$ for $c(t)=0$, which cannot be an immersion, then $\|F(x,t)\|^2 = 0$ for all $t \in [0,T)$. Which is a contradiction to eq. \eqref{ijijijijijij}. So we proved
\begin{thm}\label{nullsphere}
There are no hyperquadric homotheties of the MCF $F: M \times [0,T) \to \Rn$ with nondegenerate metric such that $F(M,0) \subset \mathcal{H}^{n-1}(0)$.
\end{thm}
\begin{flushright}
\Square
\end{flushright}
\end{subsubsection}

\begin{rem}
One could expect to find at least some stationary solutions in the light cone, like straight lines, but for such a line the metric is degenerate and thence this case is not included in Theorem \ref{nullsphere}.
\end{rem}

\begin{rem}
But, as Ecker noted in \cite{Ecker97}, the upper light-cone would immediately change to a hyperquadric and the explicitly solution to the MCF with the upper light-cone as initial condition in $\R^{1,n}$ is given by the family of graphs $\delta_t: \R^{n-1} \to \R$
$$
\delta_t (x) = \sqrt{\|x\|_{\E}^2 + 2(n-1)t},
$$
for any $t \in [0, \infty)$, which is a homothety after $t=0$.
\end{rem}

\begin{rem}
If $k>0$ then $F$ and $\H$ are pointing in oposite directions and $F_0(M)$ shrinks under the mean curvature flow.

If $k<0$ then $F$ and $\H$ are pointing in the same direction and $F_0(M)$ expands under the mean curvature flow.
%
%
\end{rem}

\begin{defn}
Let $(N,h)$ be a semi-Riemannian manifold, $M$ a smooth manifold and $F:M \to N$ an immersion such that the mean curvature vector of $F$ satisfies $\| \vec{H}(x)\|^2 \neq 0$ for all $x \in M$. The \textit{principal normal} is the vectorfield
$$
\nu := \frac{\H}{\sqrt{|\|\H\|^2|}}.
$$
\end{defn}

\begin{rem}
It is clear from equation \eqref{meancurvaturesphere} that $\H \neq 0$ everywhere for a hyperquadric homothety of the mean curvature flow and $\nabla^\perp \H = \nabla^\perp \nu = 0$.
\end{rem}

\end{subsection}

\end{section}

\begin{section}{Principal Normal Parallel in the Normal Bundle}
The two types of homotheties (self-shrinkers and self-expanders) lead, after rescaling, to different equations $\H= -F^\perp$ or $\H= F^\perp$. We restrict our attention, in this section, to the self-shrinkers of the MCF that have the principal normal parallel in the normal bundle.

If one considers a complexification of the tangent and normal bundles, $\frac{\H}{|\|\H\||}$ parallel in the normal bundle is equivalent\footnote{assuming $\|H\|^2\neq 0 \, \forall \, x \in M$} to the possibly imaginary vector field $\nu := \frac{\H}{\|\H\|}$ being parallel in the normal bundle.

In this section we prove that a compact spacelike self-shrinker cannot satisfy $\|\H\|^2 <m$ (in particular cannot be negative) for all $x \in M$ and we also prove that the being parallel in the normal bundle is enough, if the dimension of $M$ is different from 1, to ensure that a self-shrinker is hyperquadric, as the following Theorem states:

\begin{thm}\label{PrinParNorSphe}\footnote{As $\|\H\|^2 >0$ this is a slight generalization of Smoczyk's result for spacelike minimal immersed manifolds of the hyperquadrics of positive squared norm.}
Let $M$ be a closed smooth manifold and $F: M \to \Rn$ be a smooth immersion, which is a spacelike self-shrinker of the mean curvature flow, i.e. $F$ satisfies,
\begin{equation}\label{groundselfsimilar}
\hspace{15mm} \H = - F^\perp.
\end{equation}
Besides, assume $m:=\dim(M) \neq 1$. Then the mean curvature vector $\H$ satisfies  $\|\H\|^2(p) \neq 0$ for all $p \in M$ and the principal normal $\nu$ is parallel in the normal bundle ($\nabla^\perp \nu \equiv 0$) if, and only if, $F$ is a minimal immersion in the hyperquadric $\mathcal{H}^{n-1}(m)$.
\end{thm}


\begin{subsection}{Fundamental Equations}

In this subsection we calculate several equations involving the Laplacian of some tensors like the second fundamental form, the mean curvature vector, the Riemannian curvature and others. 
For this purpose we use three auxiliary tensors
$$
P_{ij}:=\langle \H, A_{ij} \rangle, \hspace{20mm} Q_{ij}:=\langle A^k_i, A_{kj} \rangle, \hspace{20mm} S_{ijkl}:=\langle A_{ij}, A_{kl} \rangle.
$$

Using Gau\ss \, equation (eq. \eqref{gauss}) we write the Ricci curvature as
\begin{equation}\label{ricciort}
R_{ij}= g^{kl}R_{kilj} = g^{kl}\langle A_{lk}, A_{ji} \rangle - g^{kl}\langle A_{jk}, A_{li} \rangle = P_{ij} - Q_{ij}.
\end{equation}

In this notation the useful Simon's equation is written as:
\begin{prop}\label{simons}
\begin{equation}
\nabla_k^\perp \nabla_l^\perp \H = \triangle^\perp A_{kl} + R_{kilj}A^{ij} - R^i_k A_{il} + Q^i_l A_{ik} - S_{kilj}A^{ij}
\end{equation}
\end{prop}
If we fix $t\in [0,T)$ the immersion $F_t$ can be constant rescaled to bring eq. \eqref{selfie} into  
$$
\H = - F^\perp.
$$

\begin{rem}
On the other hand, from Huisken (\cite{Hui90}), if an isometric immersion $G_0:M \to \Rn$ satisfies $\H_{G_0} = - G_0^\perp$ then the homothetic deformation given by
$$\label{HuiskensDeformation}
G(x,t) := \sqrt{1-2t} G_0
$$
is (up to a tangential component) the mean curvature flow, but tangential components do not change the form of the immersed manifold, so that an immersion shrinks homothetically under the MCF if\footnote{Up to rescaling}, and only if, it satisfies equation \eqref{groundselfsimilar}.
\end{rem}

We make use of the following one-form $\theta$:
\begin{equation}
 \theta := \frac{1}{2} d\|F \|^2 = \langle F_i, F \rangle dx^i
\end{equation}
such that $\theta^j F_j = \theta_i g^{ij}F_j$ is equal to $F^\top$, with $\theta_i = \langle F_i, F \rangle$. Then:
$$
\nabla_i \theta_j = \nabla_i \langle F_j, F \rangle = \langle A_{ij}, F \rangle + g_{ij}.
$$

Hence it follows
\begin{equation}
\nabla_i^\perp F^\perp = (\nabla_i(F - \theta^k F_k))^\perp = (F_i - \nabla_i\theta^k F_k - \theta^k A_{ik})^\perp = - \theta^k A_{ik}
\end{equation}
and
\begin{equation}\label{HThetaA}
\nabla_i^\perp \H = - \nabla_i^\perp F^\perp = \theta^k A_{ik}.
\end{equation}

So that
\begin{align}
\nonumber \nabla_i^\perp \nabla_j^\perp F^\perp =& - \nabla_i^\perp(\theta^k A_{jk}) = -(\nabla_i \theta^k A_{jk} + \theta^k \nabla_i A_{jk} )^\perp \\
\nonumber =& -(\nabla_i \theta^k A_{jk} + \theta^k \nabla_i A_{jk} )^\perp 
= -A_{ij} -\langle A_i^k, F^\perp \rangle A_{jk} - \theta^k \nabla_k^\perp A_{ij}
\end{align}
where we used the Codazzi equation (Theorem \ref{Codazzi Equation}) in the last step. From this follows that
\begin{equation}\label{secondderivativeH}
\nabla_i^\perp \nabla_j^\perp \H = -\nabla_i^\perp \nabla_j^\perp F^\perp = A_{ij} - P_i^k A_{kj} + \theta^k \nabla_i^\perp A_{jk}
\end{equation}
and
\begin{equation}\label{APdirectionNU}
\triangle^\perp \H = g^{ij}\nabla^\perp_i \nabla^\perp_j \H = \H - P^{ik}A_{ik} + \theta^k\nabla_k^\perp \H.
\end{equation}
Now we are able to calculate $\triangle \|\H\|^2$:
\begin{align}
\nonumber \triangle \|\H\|^2 
=& 2g^{ij}(\langle \nabla_i^\perp \nabla_j^\perp \H, \H \rangle + \langle \nabla_j^\perp \H, \nabla_i^\perp \H \rangle)\\
\nonumber =& 2\langle \triangle^\perp \H, \H \rangle + 2 \|\nabla^\perp \H\|^2 = 2 \langle \H - P^{ik}A_{ik} + \theta^k\nabla_k^\perp \H, \H \rangle + 2 \|\nabla^\perp \H\|^2\\
\label{laplaciannormmeancurvature}\triangle \|\H\|^2 =& 2\| \H \|^2 -2 \|P\|^2 + 2 \|\nabla^\perp \H\|^2 + \langle F^\top, \nabla\| \H \|^2 \rangle,
\end{align}
because
\begin{align*}
2 \langle \theta^k \nabla^\perp_k \H, \H \rangle =& 2\langle F, F_l \rangle g^{lk} \langle \nabla_k \H, \H \rangle = \langle F, F_l \rangle g^{lk} \nabla_k \langle \H, \H \rangle\\
=& \langle \langle F, F_l \rangle g^{lu} F_u, \nabla_k \langle \H, \H \rangle g^{kt} F_t \rangle = \langle F^\top, \nabla \| \H\|^2 \rangle.
\end{align*}
For $\|A\|^2$, using Simon's equation (Proposition \ref{simons}), one gets:
\begin{align}
\nonumber 2\langle A, (\nabla^\perp)^2 \H \rangle 
=& g^{tk} g^{sl} 2\langle A_{ts} ,\triangle^\perp A_{kl} + R_{kilj}A^{ij} - R^i_k A_{il} + Q^i_l A_{ik} - S_{kilj}A^{ij} \rangle\\
\nonumber =& \triangle \|A\|^2 - 2 \|\nabla^\perp A\|^2 +2R_{kilj}S^{ijkl}-2R_{ij}Q^{ij} + 2\|Q\|^2 -2S_{ikjl}S^{ijkl}.
\end{align}

On the other hand, using eq. \eqref{RicciWedge} for the Ricci tensor of the normal bundle, we get
\begin{align}
\nonumber \|R^\perp\|^2 
= \langle A_{jk}, A^j_l \rangle &\langle A^k_i, A^{li} \rangle - \langle A_{ik}, A^j_l \rangle \langle A^k_j, A^{li} \rangle -\langle A_{jk}, A^i_l \rangle \langle A^k_i, A^{lj} \rangle + \langle A_{ik}, A^i_l \rangle \langle A^k_j, A^{lj} \rangle\\
\label{ricciortogonal}\|R^\perp\|^2=& Q_{kl}Q^{kl} - S_{ikjl}S^{kjli} - S_{jkil}S^{kilj} + Q_{kl}Q^{kl} = 2\|Q\|^2 - 2S_{ikjl}S^{ijkl}.
\end{align}
So that, using these last two equations, we reach
\begin{align}
\nonumber 2\langle A, (\nabla^\perp)^2 \H \rangle =& \triangle \|A\|^2 - 2 \|\nabla^\perp A\|^2 +2R_{kilj}S^{ijkl} -2R_{ij}Q^{ij} + \|R^\perp\|^2\\
2\langle A, (\nabla^\perp)^2 \H \rangle =& \triangle \|A\|^2 - 2 \|\nabla^\perp A\|^2 + 2\|S\|^2 -2 \langle P, Q \rangle + 2 \|R^\perp\|^2,
\end{align}
where we used Gau\ss \, equation (eq. \eqref{gauss}), eq. \eqref{ricciort} and eq. \eqref{ricciortogonal} in the last step.

On the other hand, we can calculate an equation for $\triangle \|A\|^2$ using Simon's equation (Proposition \ref{simons}) in the following way: First, with eqs. \eqref{secondderivativeH} and \eqref{ricciort}, we have
\begin{align}
\nonumber \triangle^\perp A_{kl} =& \nabla_k^\perp \nabla_l^\perp \H - R_{kilj}A^{ij} + R^i_k A_{il} - Q^i_l A_{ik} + S_{kilj}A^{ij}\\
\label{triangleAkl} \triangle^\perp A_{kl} =& A_{kl} - Q_k^i A_{il} - Q^i_l A_{ik} + \theta^t \nabla_k^\perp A_{lt} + (S_{kilj}-R_{kilj})A^{ij},
\end{align}
which implies
\begin{align}
\nonumber \triangle \|A\|^2 
=& 2 \langle A_{kl} - Q_k^i A_{il} - Q^i_l A_{ik} + \theta^i \nabla_k^\perp A_{li} + (S_{kilj}-R_{kilj})A^{ij}, A^{kl} \rangle+ 2 \| \nabla^\perp A \|^2\\
\nonumber =& 2\|A\|^2 - 4 \|Q\|^2 + \langle F^\top ,\nabla \|A\|^2 \rangle + 2(2S_{kilj}-S_{klij})S^{ijkl}+ 2\| \nabla^\perp A \|^2\\
\triangle \|A\|^2 =& 2\|A\|^2 - 2 \|R^\perp\|^2 + \langle F^\top ,\nabla \|A\|^2 \rangle - 2\|S\|^2+ 2 \| \nabla^\perp A \|^2,
\end{align}
where we used Gau\ss \, equation (eq. \eqref{gauss}), equation \eqref{ricciortogonal} and
\begin{align*}
2 \langle \theta^i \nabla^\perp_k A_{li}, A^{kl} \rangle =& 2\langle F, F_t \rangle g^{ti} \langle \nabla_i A_{kl}, A^{kl} \rangle = \langle F, F_t \rangle g^{ti} \nabla_i \langle A_{kl}, A^{kl} \rangle\\
=& \langle \langle F, F_t \rangle g^{tu} F_u, \nabla_i \langle A_{kl}, A^{kl} \rangle g^{is} F_s \rangle = \langle F^\top, \nabla \| A\|^2 \rangle.
\end{align*}

\begin{thm}\label{CompactNonNegative}
Let $M$ be a closed smooth manifold and $F:M \to \Rn$ be a spacelike self-shrinker of the mean curvature flow. Then it cannot hold $\|\H\|^2<m:=dim(M)$.
\end{thm}
\begin{proof}
If $\|\H\|^2<m$ for all $x \in M$, then
\begin{equation}\label{tubdubda}
\triangle \|F\|^2 = 2 g^{ij}\langle F_i, F_j \rangle + 2\langle \triangle F, F \rangle = 2m - 2\|\H\|^2>0.
\end{equation}
But at a maximum $p$ of $\|F\|^2$ it holds $\triangle \|F\|^2 \leq 0$. Which is a contradiction.
\end{proof}

\begin{rem}
In particular there are no spacelike self-shrinkers with $\|\H\|^2<0$ and no spacelike self-shrinkers if the index of $\Rn$ is $n-m$.
\end{rem}
\end{subsection}

\begin{subsection}{The Compact Case}
Let us now consider the self-shrinkers of the MCF that satisfy the following conditions:
\begin{itemize}
 \item The mean curvature vector is not a null vector
$$
\|\H (x)\|^2 \neq 0, \mbox{ for all } x \in M.
$$
\item The principal normal $\nu:=\frac{1}{\|\H\|} \H$ is parallel in the normal bundle
$$
\nabla^\perp \nu \equiv 0,
$$
where we write $\|\H\|$ to the complex function $\sqrt{\|\H\|^2}:M \to \C$. Although Theorem \ref{CompactNonNegative} implies that $\|\H\|^2 \geq 0$ in the compact case, we also consider $\|\H\|^2 \leq 0$ as a possibility for the calculations bellow for they are of use in the non-compact case. 
\end{itemize}

The complex function $\|\H\|$ is a pure real or a pure imaginary all over $M$. So $\nu$ may not to be a real vector, but a vector field in the complexification of the pullback over $M$ of $T\R^n$ ($F^{-1}T\R^n_{\mathbb{C}}$). Over this bundle we extend the inner product and the connection linearly. 
%
%
Additionally 
for $X \in \Gamma(F^*T\R^n)$, we use $(iX)^\perp := i(X^\perp)$.

\begin{rem}
The equations considered bellow are real or pure imaginary. Thence there will be not explicit mentions of the complexifications in the calculations.
\end{rem}

A parallel principal normal (in the normal bundle) can simplify some of the previously calculated equations because of its properties:
\begin{equation}
\nabla^\perp_k \H = \nabla^\perp_k (\|\H\|\nu) = \nabla_k \|\H\| \nu
\end{equation}
and
\begin{equation}
\triangle^\perp \H = g^{ij}\nabla^\perp_i \nabla^\perp_j(\|\H\|\nu ) = g^{ij}\nabla_i \nabla_j\|\H\|\nu = \triangle \|\H\|\nu.
\end{equation}

From this, using equation \eqref{APdirectionNU}, we calculate
\begin{equation}\label{PAdirectionNU}
P^{ij}A_{ij} = \H + \theta^k \nabla^\perp_k \H - \triangle^\perp \H = (\|\H\| + \theta^k \nabla_k \|\H\| - \triangle \|\H\|) \nu,
\end{equation}
which means that $P^{ij}A_{ij}$ is in the same direction as $\nu$ (or $i \nu$, if $\nu$ is imaginary).
\begin{lem}\label{equationswithindices}
Let $F:M \to \Rn$ be an immersion such that the principal normal is parallel in normal bundle, then
$$
\begin{array}{rlrl}
 1) & P^{ij}A_{ij}= \frac{\|P\|^2}{\|\H\|}\nu & 2) & S_{ijkl}P^{ij}P^{kl}=\frac{\|P\|^4}{\|\H\|^2}\\
 3) & P_i^k A_{kj} = P_j^k A_{ki} & 4) & S_{ikjl} P^{ij} P^{kl} = Q_{il} P^i_k P^{kl}
\end{array}
$$
\end{lem}

\begin{lem}
Let $F:M \to \Rn$ be a self-shrinker of the MCF such that the principal normal is parallel in normal bundle, then
\begin{equation}
 \frac{4}{\|\H\|^4}\left< \nabla^\perp \H, \nabla^\perp A_{ij} \right> P^{ij} =\frac{2}{\|\H\|}\left< \nabla \|\H\|, \nabla\left( \frac{\|P\|^2}{\|\H\|^4} \right) \right> + 4\frac{\|P\|^2 }{\|\H\|^6} \|\nabla \|\H\|\|^2.
\end{equation}
\end{lem}
\begin{proof}
We start calculating
\begin{align*}
\left< \nabla^\perp \H, \nabla^\perp A_{ij} \right> P^{ij} =& \nabla^k \|\H\| \left< \nu, \nabla^\perp_k A_{ij} \right> P^{ij} = \nabla^k \|\H\| \nabla_k (\left< \nu, A_{ij} \right>) P^{ij}\\
=& \frac{1}{2\|\H\|}\langle \nabla\|\H\|, \nabla \|P\|^2 \rangle - \frac{\|P\|^2}{\|\H\|^2}\|\nabla \|\H\|\|^2
\end{align*}
and
\begin{align*}
\left< \nabla \|\H\|, \nabla\left( \frac{\|P\|^2}{\|\H\|^4} \right) \right> =& \left< \nabla \|\H\|, \frac{\nabla\|P\|^2}{\|\H\|^4} - \frac{4\|P\|^2 \|\H\|^3 \nabla\|\H\|}{\|\H\|^8} \right>\\
=& \frac{\left< \nabla \|\H\|, \nabla\|P\|^2 \right>}{\|\H\|^4} - 4\frac{\|P\|^2 }{\|\H\|^5} \|\nabla \|\H\|\|^2.
\end{align*}
These two equations imply that
\begin{equation}
\label{nablaHnablaAij}\frac{4}{\|\H\|^4}\left< \nabla^\perp \H, \nabla^\perp A_{ij} \right> P^{ij}= \frac{2}{\|\H\|}\left< \nabla \|\H\|, \nabla\left( \frac{\|P\|^2}{\|\H\|^4} \right) \right> + 4\frac{\|P\|^2 }{\|\H\|^6} \|\nabla \|\H\|\|^2.
\end{equation}
\end{proof}
We continue by calculating $\diamond:=\frac{2}{\|\H\|^4} \left\| \nabla_i\|\H\|\frac{P_{jk}}{\|\H\|} - \|\H\| \nabla_i\left(\frac{P_{jk}}{\|\H\|}\right) \right\|^2$,
\begin{align}
\nonumber \diamond =& \frac{2}{\|\H\|^6}\|\nabla\|\H\|\|^2\|P\|^2 + \frac{2}{\|\H\|^2}\left\|\nabla_i\left(\frac{P_{jk}}{\|\H\|}\right) \right\|^2 - \frac{4}{\|\H\|^4}\nabla_i \|\H\| \nabla^i\left(\frac{P_{jk}}{\|\H\|}\right) P^{jk}\\
%
\nonumber =& \frac{2}{\|\H\|^2}\left\|\nabla_i\left(\frac{P_{jk}}{\|\H\|}\right) \right\|^2 - 2\frac{\|P\|^2}{\|\H\|^6}\|\nabla\|\H\|\|^2 - \frac{2}{\|\H\|}\left<\nabla \|\H\|, \nabla \left(\frac{\|P\|^2}{\|\H\|^4} \right)\right>
\end{align}
and
\begin{align*}
\frac{2}{\|\H\|^2} \left\|\nabla \left(\frac{P}{\|\H\|} \right) \right\|^2=& \frac{2}{\|\H\|^2}\left\|\frac{\nabla_i P_{jk}}{\|\H\|} - \frac{\nabla_i \|\H\| P_{jk}}{\|\H\|^2}\right\|^2\\
=&2\frac{\|\nabla P\|^2}{\|\H\|^4} -6 \frac{\|P\|^2}{\|\H\|^6}\|\nabla\|\H\|\|^2 - \frac{2}{\|\H\|}\left<\nabla \|\H\|, \nabla\left(\frac{\|P\|^2}{\|\H\|^4}\right) \right>.
\end{align*}
With this we get equation:
\begin{align}
\label{nablaHandP} \frac{2}{\|\H\|^4} \left\| \nabla_i\|\H\|\frac{P_{jk}}{\|\H\|} - \|\H\| \nabla_i\left(\frac{P_{jk}}{\|\H\|}\right) \right\|^2=& 2\frac{\|\nabla P\|^2}{\|\H\|^4} - 8\frac{\|P\|^2}{\|\H\|^6}\|\nabla\|\H\|\|^2\\
\nonumber &- \frac{4}{\|\H\|}\left<\nabla \|\H\|, \nabla \left(\frac{\|P\|^2}{\|\H\|^4} \right)\right>.
\end{align}

\begin{lem}\label{laplacianQuocient}
Let $F:M \times [0,T) \to \Rn$ be a self-shrinker of the MCF such that $\|\H\|^2 \neq 0$ for all $x\in M$ and the principal normal is parallel in the normal bundle. Then
\begin{align}
\triangle \left(\frac{\|P\|^2}{\|\H\|^4} \right) =& \frac{2}{\|\H\|^4} \left\|\nabla_i\|\H\|\frac{P_{jk}}{\|\H\|} - \|\H\| \nabla_i \left(\frac{P_{jk}}{\|\H\|} \right) \right\|^2 \\
\nonumber &+ \left<F^\top, \nabla\left(\frac{\|P\|^2}{\|\H\|^4} \right) \right> - \frac{2}{\|\H\|}\left<\nabla\|\H\|, \nabla\left(\frac{\|P\|^2}{\|\H\|^4} \right)\right>.
\end{align}
\end{lem}
\begin{proof}
We begin using equations \eqref{APdirectionNU} and \eqref{triangleAkl} to calculate
\begin{align*}
\triangle P_{ij} =& \nabla^k \nabla_k \langle \H, A_{ij} \rangle = \langle {\nabla^k}^\perp \nabla_k^\perp\H, A_{ij} \rangle + 2 \langle \nabla_k^\perp \H, {\nabla^k}^\perp A_{ij} \rangle + \langle \H, {\nabla^k}^\perp \nabla_k^\perp A_{ij} \rangle\\
=& 2 P_{ij} +2 \langle \nabla^\perp \H, \nabla^\perp A_{ij}\rangle - Q_i^k P_{kj}- Q_j^k P_{ki} + 2(S_{ikjl}- S_{ijkl})P^{kl} + \langle F^\top, \nabla P_{ij} \rangle,
\end{align*}
where we used Gau\ss \, equation (eq. \eqref{gauss}) and $\theta^k \nabla_k P_{ij} = \langle F^\top, \nabla P_{ij} \rangle$. From this follows that
\begin{align*}
\triangle \|P\|^2 =& \triangle(P_{ij}P^{ij}) = 2\triangle P_{ij} P^{ij} +2\langle\nabla P, \nabla P \rangle\\
=& 2 \|\nabla P\|^2 + \langle F^\top, \nabla \|P\|^2\rangle + 4 \langle \nabla^\perp \H, \nabla^\perp A_{ij}\rangle P^{ij}- 4 \frac{\|P\|^4}{\|\H\|^2} + 4 \|P\|^2,
\end{align*}
where we used Lemma \ref{equationswithindices}.

On the other side
\begin{align*}
\triangle \left(\frac{\|P\|^2}{\|\H\|^4}\right) 
=& \frac{\triangle \|P\|^2}{\|\H\|^4} - 8\frac{\nabla^i \|\H\|}{\|\H\|}\left(\frac{\nabla_i\|P\|^2}{\|\H\|^4}- 4\frac{\|P\|^2\nabla_i\|\H\|}{\|\H\|^5}\right)\\
&-2\frac{\|P\|^2}{\|\H\|^6}\left(2\triangle \|\H\| \|\H\| + 2\|\nabla\|\H\|\|^2\right) -8\frac{\|P\|^2 \|\nabla \|\H\|\|^2}{\|\H\|^6},
\end{align*}
which implies that
$$
\triangle \left(\frac{\|P\|^2}{\|\H\|^4}\right) = \frac{\triangle \|P\|^2}{\|\H\|^4} - \frac{8}{\|\H\|}\left<\nabla\|\H\|, \nabla \frac{\|P\|^2}{\|\H\|^4}\right> - 2\frac{\|P\|^2}{\|\H\|^6}\triangle \|\H\|^2-8\frac{\|P\|^2 \|\nabla \|\H\|\|^2}{\|\H\|^6}.
$$

Using the equations for $\triangle \|P\|^2$ and $\triangle \|\H\|^2$ (eq. \eqref{laplaciannormmeancurvature}) we get,
\begin{align*}
\triangle \left(\frac{\|P\|^2}{\|\H\|^4}\right) 
=& \left<F^\top, \nabla \left(\frac{\|P\|^2}{\|\H\|^4} \right) \right> + \frac{2 \|\nabla P\|^2 + 4 \langle \nabla^\perp \H, \nabla^\perp A_{ij}\rangle P^{ij}}{\|\H\|^4}\\
&- \frac{8}{\|\H\|}\left<\nabla\|\H\|, \nabla \left(\frac{\|P\|^2}{\|\H\|^4} \right)\right> - 12\frac{\|P\|^2}{\|\H\|^6} \|\nabla \|\H\|\|^2,
\end{align*}
then we apply eqs. \eqref{nablaHandP} and \eqref{nablaHnablaAij} together with $\nabla^\perp \H = \nabla \|\H\| \nu$ to prove the Lemma.
\end{proof}

\begin{prop}
Let $M$ be a closed smooth manifold and $F: M \to \Rn$ be a smooth immersion, which is a spacelike self-shrinker of the mean curvature flow. Besides, assume that the mean curvature vector $\H$ satisfies $\|\H\|^2 \neq 0$ and the principal normal $\nu$ satisfies $\nabla^\perp \nu =0$. Then
\begin{equation}\label{result}
\left\| \nabla_i\|\H\|\frac{P_{jk}}{\|\H\|} - \|\H\| \nabla_i\left(\frac{P_{jk}}{\|\H\|}\right)\right\|^2 =0.
\end{equation}
\end{prop}
\begin{proof}
Although the function $\|\H\|$ may be imaginary, the 3-tensor in eq. \eqref{result} is real. Then
$$
\left\| \nabla_i\|\H\|\frac{P_{jk}}{\|\H\|} - \|\H\| \nabla_i\left(\frac{P_{jk}}{\|\H\|}\right)\right\|^2 \geq 0.
$$

From Lemma \ref{laplacianQuocient}, we can write 
$$
\triangle \left(\frac{\|P\|^2}{\|\H\|^4} \right) \leq \left<F^\top - \frac{2}{\|\H\|}\nabla\|\H\|, \nabla\left(\frac{\|P\|^2}{\|\H\|^4} \right) \right>.
$$
The strong elliptic maximum principle implies that $u$ is constant.
Then $\nabla \left(\frac{\|P\|^2}{\|\H\|^4}\right) =0$ and $\triangle \left(\frac{\|P\|^2}{\|\H\|^4}\right) =0$. Hence theorem \ref{laplacianQuocient} implies that
$$
\left\| \nabla_i\|\H\|\frac{P_{jk}}{\|\H\|} - \|\H\| \nabla_i\left(\frac{P_{jk}}{\|\H\|}\right)\right\|^2 =0.
$$
\end{proof}

We now rewrite the equality that we just proved in another way:

First, eq. \eqref{result} implies
\begin{equation}\label{locozero}
\nabla_i\|\H\|\frac{P_{jk}}{\|\H\|} - \|\H\| \nabla_i\left(\frac{P_{jk}}{\|\H\|}\right) =0,
\end{equation}
as a tensor, because this is a covariant tensor over $M$ and $M$ is spacelike.

Second, using the Codazzi equation (eq. \eqref{Codazzi Equation}) and $\nabla^\perp \nu = 0$, we calculate
$$
\nabla_i \left(\frac{P_{jk}}{\|\H\|} \right) = \nabla_i \langle \nu, A_{jk} \rangle = \langle \nu, \nabla_i^\perp A_{jk} \rangle = \langle \nu, \nabla_j^\perp A_{ik} \rangle = \nabla_j \langle \nu, A_{ik} \rangle = \nabla_j \left(\frac{P_{ik}}{\|\H\|} \right).
$$

Third, using equation \eqref{locozero}, we write
\begin{align*}
\nabla_i\|\H\|\frac{P_{jk}}{\|\H\|} - \|\H\| \nabla_i\left(\frac{P_{jk}}{\|\H\|}\right) 
=& \nabla_i\|\H\|\frac{P_{jk}}{\|\H\|} - \nabla_j \|\H\| \frac{P_{ik}}{\|\H\|},
\end{align*}
which implies
$$
0 = \left\| \nabla_i\|\H\|\frac{P_{jk}}{\|\H\|} - \|\H\| \nabla_i\left(\frac{P_{jk}}{\|\H\|}\right)\right\|^2 = \left\| \nabla_i\|\H\|\frac{P_{jk}}{\|\H\|} - \nabla_j \|\H\| \frac{P_{ik}}{\|\H\|}\right\|^2.
$$

Now, expanding this norm we find
\begin{align}
\label{solutionspherical} \|\nabla \|\H\| \|^2 \|P\|^2 - \| \nabla_i \|\H\| P_k^i \|^2 =& 0.
\end{align}
With this formula we can show that $F$ is hyperquadric, i. e. $\|F\|^2 = q \in \R$.

\noindent\textbf{What remains to prove of Theorem \ref{PrinParNorSphe}} \textit{Let $M$ be a closed smooth manifold and $F: M \to \Rn$ be an immersion, which is a spacelike self-shrinker of the mean curvature flow. Besides, assume that the mean curvature vector $\H$ satisfies $\|\H\|^2 \neq 0$ and the principal normal $\nu$ satisfies $\nabla^\perp \nu =0$. If $m:=\dim(M) \neq 1$, then}
$$
\|F(x)\|^2 = m\, \forall \,x \,\in \, M.
$$
\begin{proof}
We now calculate at a point $p \in M$ fixed. As the 2-tensor $P$ is symmetric, it is also diagonalizable and has only real eigenvalues $\lambda_1, \ldots, \lambda_m$. Let $V_1, \ldots, V_m$ be an orthonormal basis of eigenvectors associated with $\lambda_1, \ldots, \lambda_m$. Then we write $\nabla \|\H\| = \sum_i \alpha_i V_i$, $\alpha_i \in \C$ so that by equation \eqref{solutionspherical}
\begin{equation}\label{patapatapatapon}
0 = \|P\|^2 \|\nabla \|\H\|\|^2 - \|P(\nabla \|\H\|)\|^2 = \sum_i \lambda_i^2 (\|\nabla \|\H\|\|^2 - \alpha_i^2),
\end{equation}
but $\lambda_i^2 \geq 0$ because $\lambda_i \in \R$, beyond this $\|\H\|$ is pure real or pure imaginary everywhere and all the $\alpha_i$'s have to agree with $\|\H\|$ about being real or imaginary, which implies
$$
\|\nabla \|\H\|\|^2 - \alpha_i^2 = \sum_{j \neq i} \alpha^2_j
$$
being nonnegative for all $i \in \{1, \ldots, m\}$ if $\|\H\|$ is real or nonpositive if $\|\H\|$ is imaginary. This implies, with eq. \eqref{patapatapatapon}, that
$$
\lambda_i^2 (\|\nabla \|\H\|\|^2 - \alpha_i^2) = 0 \,\,\, \forall i \in \{1, \ldots, m\}.
$$
As $tr (P) = P_{ij}g^{ij} = \|\H\|^2 \neq 0$, it follows that $P\neq 0$ and there is at least one $j \in \{1, \ldots, n\}$ such that $\lambda_j \neq 0$ and the last equation shows that
$$
0 = \|\nabla \|\H\|\|^2 - \alpha_j^2 = \sum_i \alpha_i^2 - \alpha_j^2 = \sum_{i \neq j} \alpha_i^2 \Longrightarrow \alpha_i = 0 \,\, \forall i \neq j,
$$
because the $\alpha_i$'s are all real or all imaginary. From this follows that $\| \nabla \|\H\|\|^2 = \alpha_j^2$ and $\nabla \|\H\| = \alpha_j V_j$.

Assume that there is an $x \in M$ such that $\nabla \|\H\|(x) \neq 0$.

Then $\alpha_j \neq 0$ and for all $i \neq j$
$$
0= \lambda_i^2 (\| \nabla \|\H\| \|^2 -\alpha_i^2) = \lambda_i^2 \alpha_j^2 \Longrightarrow \lambda_i = 0,
$$
so $P_{ij}$ has only one nonzero eigenvalue and the associated eigenvector is $\nabla \|\H\|/\|\nabla \|\H\|\|$.

At this point we have
$$
\|P\|^2 = \lambda_j^2 = (\tr P)^2 = \|\H\|^4 \Longrightarrow \frac{\|P\|^2}{\|\H\|^4} = 1,
$$
but we have already shown that this quotient is constant, so that the equation $\frac{\|P\|^2}{\|\H\|^4} = 1$ holds not only at this point but everywhere in $M$.

Then, using $\|P\|^2=\|\H\|^4$, with equation \eqref{laplaciannormmeancurvature} we calculate
\begin{align*}
 2\|\H\| \triangle \|\H\| =& 2\| \H \|^2 -2 \|\H\|^4 + 2 \|\H\|\langle F^\top, \nabla \| \H \| \rangle
\end{align*}
and it follows
\begin{equation}
\triangle \|\H\| = \| \H \| - \|\H\|^3 + \langle F^\top, \nabla \|\H \| \rangle.
\end{equation}

We integrate both sides of this equation. First integrate the terms of it separately taking advantage of the fact that $M$ is closed:
$$
\int_M \triangle \|\H\| = 0,
$$
because of the Divergence Theorem
, and
\begin{align*}
\int_M \langle F^\top, \nabla\|\H\| \rangle =& \int_M \langle F, F_l \rangle g^{lk} \nabla_k \|\H\| = -\int_M \nabla_k \langle F, F_l \rangle g^{lk} \|\H\| \\
=& -m \int_M \|\H\| + \int_M \|\H\|^3,
\end{align*}
such that
$$
0 = \int_M \triangle \|\H\| = \int_M \| \H \| - \|\H\|^3 + \langle F^\perp, \nabla |\H | \rangle = (1-m)\int_M \|\H\|,
$$
which is impossible for $m \neq 1$.

From this contradiction we know that $\nabla \|\H\| = 0$ everywhere in $M$. It follows that $\nabla^\perp \H = \nabla \|\H\| \nu = 0$ and that the norm of $\H$ is constant.

On the other hand
\begin{equation}\label{terminalogo}
\triangle \|F\|^2 = 2\langle F_i, F_j \rangle g^{ij} + 2\langle \triangle F, F \rangle = 2g_{ij} g^{ij} + 2\langle \H, F \rangle = 2m -2\|H\|^2.
\end{equation}

If the constant $2m -2\|H\|^2$ is other than zero (for example $>0$) it would lead to a contradiction with the second derivative's test 
, so that $\triangle \|F\|^2 = 0$ everywhere in $M$.

Again using the maximum principle, we find that $\|F\|^2$ is constant. This norm can be calculated seeing that $\left<F_i, F\right> = 0$, which implies that $F \in \Gamma(TM^\perp)$, so that $\H = -F$ and replacing $\|F\|^2 = \|\H\|^2$ and $\triangle \|F\|^2 = 0$ in eq. \eqref{terminalogo} we get $\|F\|^2=\|\H\|^2=m$.
\end{proof}
Note that the condition dim$(M) \neq 1$ is optimal, because the result does not hold for the curve shortening flow, then the Abresch \& Langer curves are not contained in a circle.

%
\end{subsection}

\end{section}

\begin{section}{The Non-Compact Case}
We now consider non-compact self-shrinkers and need to integrate over $M$ with respect to a backwards heat kernel. 
Let us consider in $\R^n$ the usual topology. A set $B \subset \R^n$ is \textit{unbounded} if there is no compact set containing $B$.


\begin{rem}
In the pseudo-euclidean case there are minimal submanifolds of the hyperquadrics, which are noncompact and are homotheties of the mean curvature flow with principal normal parallel in the normal bundle. These hyperquadrics are asymptotic to the light cone and, in particular, have the norm $\|F\|^2$ bounded, thence they do not satisfy the conditions needed to integrate and do not appear in our results.
\end{rem}

In the compact case we proved that $\|\H\|^2 < m$ implies that $F$ is not a self-shrinker of the MCF; in the non-compact case a similar result holds.

\begin{thm}\label{NONCompactNonNegative}
The mean curvature vector of a stochastic complete, spacelike, self-shrinker of the mean curvature flow $F:M \to \Rn$ cannot satisfy, for all $p \in M$,
$$
\|\H\|^2 < m - \epsilon,
$$
for some $\epsilon >0$ if $\sup_M \|F\|^2< +\infty$.
\end{thm}
\begin{proof}
If there is an $\epsilon > 0$ such that $\|\H\|^2<m - \epsilon$ for all $x \in M$, then
$$
\triangle \|F\|^2 = 2 g^{ij}\langle F_i, F_j \rangle + 2\langle \triangle F, F \rangle = 2m - 2\|\H\|^2 > 2\epsilon,
$$
but by the weak Omori-Yau maximum principle there is a sequence $\{x_k\} \subset M$ with
$$
\triangle \|F\|^2(x_k) \leq \frac{1}{k},
$$
which contradicts $\triangle \|F\|^2(x) > 2\epsilon$ for all $x \in M$.
\end{proof}

\begin{rem}
In particular, there are no stochastic complete, spacelike self-shrinkers of the mean curvature flow with $\sup_M \|F\|^2< +\infty$ and $\|\H\|^2 \leq 0.$
\end{rem}


\begin{defn}
Let $\Rn$ be an inner product space and $\{e_1, \ldots, e_n\}$ an orthonormal basis such that $\langle e_\alpha, e_\alpha \rangle = -1$ for $\alpha \in \{1, \ldots, q\}$ and $\langle e_\alpha, e_\alpha \rangle = 1$ for $\alpha \in \{q+1, \ldots, n\}$, which we denote $\R^{q,n}$. For a vector $X \in \R^{q,n}$ we define ($X_-$) and ($X_+$) as the projections of $X$ in span$\{e_1, \ldots, e_q\}$ and in span$\{e_{q+1}, \ldots, e_n\}$ respectively.
\end{defn}

\begin{defn}\label{mainly}
Let $M$ be a smooth manifold and $F:M \to \R^{q,n}$ be an immersion with $F(M)$ unbounded. We say that $F$ (or $F(M)$) is \textit{mainly positive} if there is an $\epsilon >0$ and $k \in \R$, such that $\forall x \in M$
$$
\|F(x)\|^2_\E \geq k \Longrightarrow  -\frac{\|F(x)_-\|^2}{\|F(x)_+\|^2} \leq 1 - \epsilon.
$$
And we say that $F$ (or $F(M)$) is \textit{mainly negative} if there is an $\epsilon >0$ and $k \in \R$, such that $\forall x \in M : \|F(x)\|^2_\E \geq k \Longrightarrow  -\frac{\|F(x)_+\|^2}{\|F(x)_-\|^2} \leq 1 - \epsilon.$
\end{defn}

\begin{multicols}{2}
\begin{minipage}{\linewidth}
\centering\includegraphics[width=42mm, angle=90]{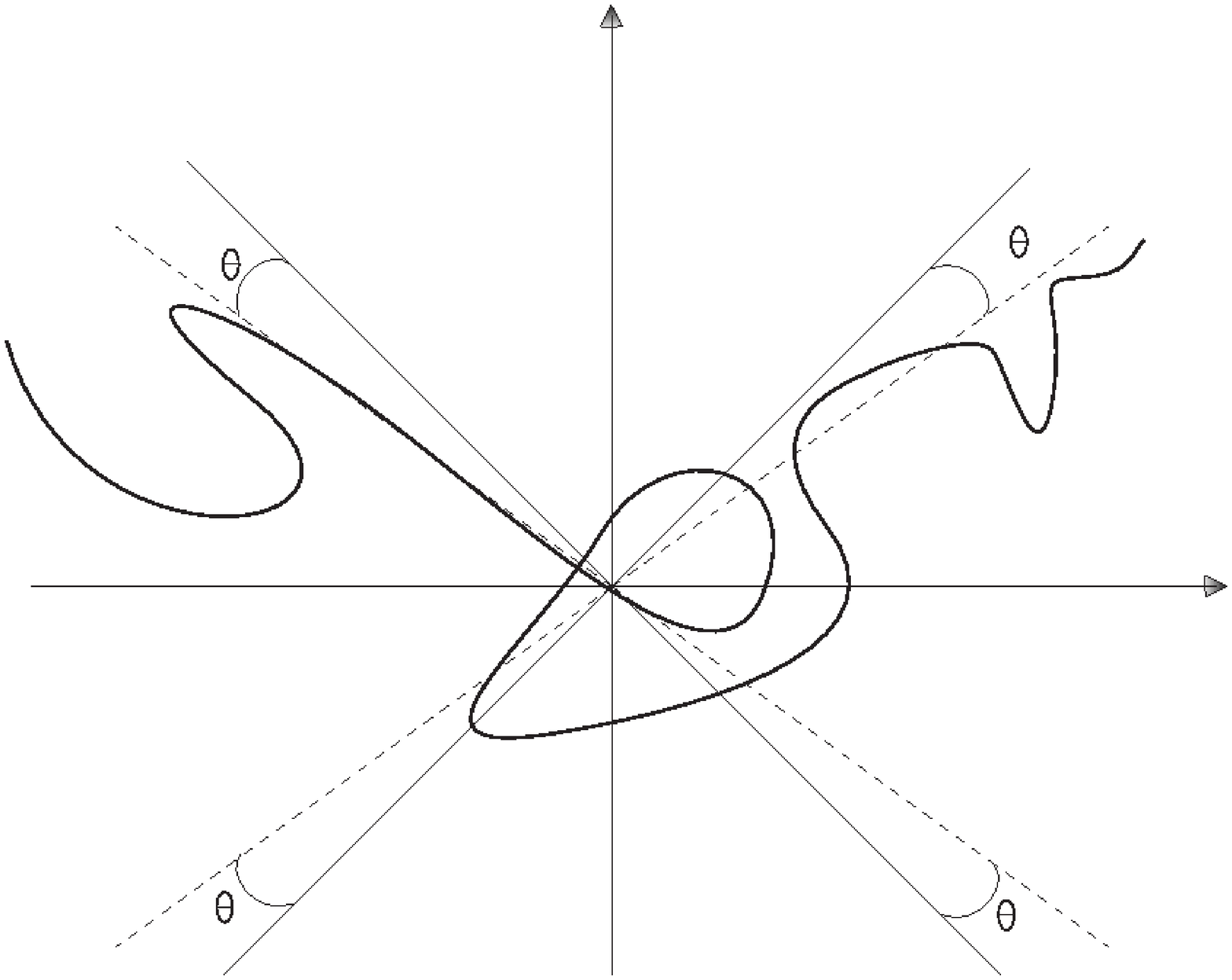}
    \end{minipage}
    
This means that there is an (Euclidean) angle $\theta$ with $\tan \left(\frac{\pi}{4}-\theta\right) < 1-\epsilon$ between $F(x)$ and the light cone for any $x\in M$ such that $F(x)$ lies outside some big euclidean sphere (or $\tan \left(\frac{\pi}{4}-\theta\right) > 1 + \epsilon$ in the mainly negative case).
\end{multicols}


\begin{lem}\label{chuchuchupiuiiiii}
If $F(M)$ is mainly positive and unbounded, then $\|F\|^2 \geq \frac{\epsilon}{2}\|F(x)\|^2_\E $ and $\|F\|^2$ is unbounded.
\end{lem}

\begin{lem}\label{iwillfindyou}
If $F(M)$ is mainly negative and unbounded, then $-\|F\|^2 \geq \frac{\epsilon}{2}\|F(x)\|^2_\E $ and $\|F\|^2$ is unbounded.
\end{lem}

Now we consider the behavior of $\|F(x)\|^2$ for $x$ in $M$ satisfying def. \ref{mainly}.

\begin{rem}
If $F$ is a spacelike self-shrinker such that $F(M)$ is mainly negative and unbounded, then for $x \in M$ with $\|F(x)\|^2_\E > k$, for $k$ as in def. \ref{mainly}, it holds that
$$
0 > \|F(x)\|^2 = \|F^\perp(x)\|^2 + \|F^\top(x)\|^2 \geq \|\H(x)\|^2,
$$
but if $M$ is stochastic complete, then Theorem \ref{NONCompactNonNegative} implies that $F$ cannot be a self-shrinker of the MCF with $\|\H\|^2 (p) \neq 0$ for all $p\in M$.
\end{rem}

%
In order to integrate we need further assumptions on $F$:

\begin{defn}\label{tuningo}
Let $F:M \to \R^{q,n}$ be a spacelike isometric immersion. We say that $F$ has \textit{bounded geometry} if:
\begin{enumerate}
\item There are $c_k, d_k \in \R$ for every $k \in \mathbb{N} \cup \{0\}$ such that
\begin{align*}
\|(\nabla)^k A_+\|^2 \leq c_k,\\
-\|(\nabla)^k A_-\|^2 \leq d_k.
\end{align*}
\item The function $\frac{1}{\|\H\|}$ grows polynomially with respect to $\|F\|^2$.
\item The growth of volume of geodesic balls and their boundaries is polynomial with respect to the radius.
\item $F$ is \textit{inverse Lipschitz} with respect to the euclidean norm in $\R^n$.
\end{enumerate}
\end{defn}

Then the bounded geometry assumption \textbf{excludes the mainly negative case} for spacelike self-shrinkers, because $-d_0 \leq \|\H\|^2 \leq c_0$ and, in this case, $\|F\|^2$ has no lower bound (by Lemma \ref{iwillfindyou} ), but $\|F\|^2 = \|F^\top\|^2 + \|F^\perp\|^2$ and $\|F^\top\|^2 \geq 0$, which implies that $\|F^\perp\|^2$ is not bounded below. This contradicts $\|F^\perp\|^2 = \|\H\|^2$. So that:

\begin{thm}\label{noMainlyNegative}
There are no unbounded mainly negative spacelike self-shrinkers of the mean curvature flow with bounded geometry.
\end{thm}

Some control on geodesic balls of $M$ is needed.

\begin{lem}\label{LemmaDoido}
Let $F:M \to \R^{q,n}$ be an inverse Lipschitz immersion with respect to the euclidean norm in $\R^n$, $\Omega_R := \{X \in F(M) \subset \R^{q,n} : \|X\|_\E < R\}$ and $p \in M$ be a fixed point such that $F(p) \in \Omega_R$. Then there is a geodesic ball $B_{R'}(p)$ of radius $R'=2R/k$, where $k$ is the constant in the inverse Lipschitz condition, such that $F^{-1}(\Omega_R) \subset B_{R'}(p)$.
\end{lem}

\begin{cor}\label{narutouzumaki}
Let $\Omega$, $R$,  $p \in F^{-1}(\Omega)$ and $R'$ be as in the last Lemma and $y \in M$. Then
$$
d(p,y)> R' \hspace{5mm} \Rightarrow \hspace{5mm} y \notin \Omega,
$$
this means $\|F(y)\|_\E > R$.
\end{cor}

We continue by proving results for mainly positive immersions.

\begin{rem}\label{seilaoque}
From Lemma \ref{chuchuchupiuiiiii}: $\|F(x)\|^2_\E \leq \frac{2}{\epsilon}\|F(x)\|^2$.
\end{rem}

We get a polynomial control of the radius of big geodesic balls in terms of $\|F\|^2$:

\begin{lem}\label{taxingo}
Let $F:M \to \R^{q,n}$ be a mainly positive, inverse Lipschitz immersion and $p \in M$. Then there are $R, k_1, k_2 \in \R$ such that $x \notin B_p(R) \Rightarrow d(p,x) \leq k_1\|F(x)\| + k_2$.
\end{lem}

Although it is necessary that $\|F\|^2 \to +\infty$ note that $\|\H\|^2$ could still be negative.



\begin{lem}\label{taracotaco}
For any $X,Y \in \R^{q,n}$ it holds
$$
| \langle X, Y \rangle | \leq \|X_+\| \|Y_+\| + \sqrt{(-\|X_-\|^2)} \sqrt{(-\|Y_-\|^2)}
$$
\end{lem}

This implies:

\begin{lem}\label{trecotreco}
If $A,B \in \Gamma(F^*T\R^{q,n} \otimes TM \otimes \ldots \otimes TM \otimes T^*M \otimes \ldots \otimes T^*M)$ are such that $\|A_+\|, \|A_-\|, \|B_+\|, \|B_-\|$ grow polynomially with $\|F\|^2$, then so does $|\langle A, B \rangle|$.
\end{lem}

\begin{rem}
At any point $p \in M$, 
\begin{equation}\label{oioioioioioioioi}
\|F\|^2 = \|\H\|^2 +\|F^\top\|^2,
\end{equation}
so that, from $|\|\H\|^2|\leq c_0 +d_0$, it holds that $\|F^\top\|^2$ grows polynomially with $\|F\|^2$.
\end{rem}

\begin{lem}\label{ninjajiraia}
Let $F:M \to \R^{q,n}$ be a spacelike, mainly positive, immersion with bounded geometry and $f:M \to \R$ be some kind of polynomial (of inner products) of $\H$, $A$, their covariant derivatives, $F$, $F^\top$ and the function $\frac{1}{\|\H\|}$, then $f$ has polynomial growth with respect to $\|F\|^2$.
\end{lem}

We will integrate over the whole manifold with respect to the following \textit{heat kernel}: $\rho: M \to \R$ defined as
$$
\rho(x) := \exp\left(-\frac{\|F\|^2}{2} \right).
$$

\begin{lem}\label{nasultimas}
Let $F:M \to \R^{q,n}$ be a spacelike, mainly positive, immersion with bounded geometry and $F(M)$ unbounded, beyond this let $f:M \to \R$ be some polynomial (of inner products) of $\H$, $A$, their covariant derivatives, $F$, $F^\top$ and the function $\frac{1}{\|\H\|}$. Then
$$
\left|\int_M f \rho d\mu\right| < \infty;
$$
beyond this, one can use partial integration
$$
\int_M \rho \mbox{div} (\nabla f(x)) d\mu =- \int_M \langle \nabla \rho, \nabla f(x)\rangle d\mu.
$$
\end{lem}

By this, all integrals in the next Lemma are finite.

\begin{lem}\label{lemaImportante}
Let $F:M \to \R^{q,n}$ be a spacelike, mainly positive, self-shrinker of the mean curvature flow with bounded geometry such that $F(M)$ is unbounded. Beyond this, let $F$ satisfy $\|\H\|^2 \neq 0$ and $\nabla^\perp \nu =0$. Then
$$
\nabla_i\|\H\|\frac{P_{jk}}{\|\H\|} - \|\H\| \nabla_i \left(\frac{P_{jk}}{\|\H\|} \right) =0.
$$
\end{lem}

\begin{proof}
The expression
$$
\int_M \rho \frac{\|P\|^2}{\|\H\|^2} \triangle \left(\frac{\|P\|^2}{\|\H\|^4} \right) d\mu
$$
can be calculated using partial integration or Lemma \ref{laplacianQuocient}. Equaling these two one finds
\begin{equation}\label{importantintegral}
\int_M 2\rho \frac{\|P\|^2}{\|\H\|^6} \left\|\nabla_i\|\H\|\frac{P_{jk}}{\|\H\|} - \|\H\| \nabla_i \left(\frac{P_{jk}}{\|\H\|} \right) \right\|^2 +\rho \|\H\|^2 \left\|\nabla \left(\frac{\|P\|^2}{\|\H\|^4} \right)\right\|^2 d\mu =0.
\end{equation}
but the two summands inside the integral have the same sign everywhere. This implies in particular, using $\|P\|^2 \neq 0$ (because $P =0$ would imply $\|\H\|=0$), that
$$
\nabla_i\|\H\|\frac{P_{jk}}{\|\H\|} - \|\H\| \nabla_i \left(\frac{P_{jk}}{\|\H\|} \right) =0.
$$
\end{proof}

And we have the same result as eq. \ref{result} in the compact case. Then we follow exactly as in the previous section (the compact case) to get:

\begin{lem}\label{cucuxico}
Let $M$ be a smooth manifold and $F:M \to \R^{q,n}$ be mainly positive, spacelike, self-shrinker of the mean curvature flow with bounded geometry such that $F(M)$ is unbounded, beyond that let $F$ satisfy the conditions: $\|\H\|^2(p) \neq 0$ for all $p\in M$ and the principal normal is parallel in the normal bundle ($\nabla^\perp \nu \equiv 0$). Then one of the two holds
\begin{enumerate}
 \item $\nabla \|\H\| = 0$ everywhere on $M$\\
 \item \vspace{-3mm}There is a point $p \in M$ with $\nabla \|\H\|(p) \neq 0$, at which $\frac{\nabla \|\H\|}{\|\nabla \|\H\|\|}$ is the only eigenvector associated with a nonzero eigenvalue of $P$.
\end{enumerate}
\end{lem}

We have to treat these two cases separately.

\begin{section}{The First Case}



\begin{thm}\label{ResultNoncompact1}
Let $M$ be a smooth manifold and $F:M \to \R^{q,n}$ be a mainly positive, spacelike, shrinking self-similar solution of the mean curvature flow with bounded geometry such that $F(M)$ is unbounded. Beyond that, let $F$ satisfy the conditions: $\|\H\|^2(p) \neq 0$, $\forall p \in M$, and the principal normal is parallel in the normal bundle ($\nabla^\perp \nu \equiv 0$). If $\nabla \|\H\|(p) = 0$ for all $p \in M$, then
\begin{equation}\label{result1}
F(M) = \mathcal{H}_r \times \R^{m-r},
\end{equation}
where $\mathcal{H}_r$ is an $r$-dimensional minimal surface of the hyperquadric $\mathcal{H}^{n-1}(r)$ with $\|\H\|^2=r>0$ and $R^{m-r}$ is an $m-r$ dimensional spacelike affine space in $\R^{q,n}$.
\end{thm}

\begin{proof}
First we see that $\nabla \|\H\| = 0$ implies $\nabla^\perp \H = \nabla \|\H\|\nu = 0$ and, with eq. \eqref{HThetaA},
\begin{equation}\label{thetaA}
\theta^i A_{ij}=0.
\end{equation}
On the other hand, $\nabla \|\H\|=0$ implies that $\|\H\|^2$ is constant, so that, with Lemma \ref{lemaImportante}, it holds $\nabla P = 0$ and then equation \eqref{secondderivativeH} implies
\begin{align}
\nonumber \langle \nabla_i^\perp \nabla_j^\perp \H, \H \rangle &= \langle A_{ij} - P_i^k A_{kj} + \theta^k \nabla_i^\perp A_{jk}, \H \rangle\\
\label{P,P2}0 &= P_{ij} - P_i^kP_{kj},
\end{align}
so that $P = P^2$; i. e. $P$ is a projection and can only have 1 and 0 as eigenvalues.

Because of $\nabla_k P_{ij} = 0$ we get $\nabla_k \|P\|^2(p) = 0$, but $\|P\|^2(p)$ is equal to the number of eigenvalues 1, thus their number is constant and 
\begin{equation}\label{normHsquaredequalr}
\|\H\|^2 = \mbox{tr}P = r>0.
\end{equation}

We consider the eigenspaces associated with these two eigenvalues, they define the distributions $\mathcal{E}M$ and $\F M$ given, at any point $p\in M$, by
\begin{equation}
\mathcal{E}M_p := \{V\in T_pM : P(V) = V\}, \hspace{5mm} \mathcal{F}M_p := \{V\in T_pM : P(V) = 0\},
\end{equation}
or in local coordinates $P^j_i V^i = V^j$ (and $P^j_i V^i = 0$). As the eigenspaces are orthogonal we have $T_pM = \mathcal{E}M_p \oplus \mathcal{F}M_p$. From equation \eqref{thetaA} we have, for all $V \in\mathcal{E}_pM$, that 
\begin{equation}\label{kernelTheta}
\theta(V) = \theta_j V^j = \theta_j P^j_iV^i =0,
\end{equation}
which means that $\mathcal{E}M_p \subset \ker(\theta)$.


\begin{lem}
Under the hypothesis of Theorem \ref{ResultNoncompact1} the distributions $\mathcal{E} M$ and $\F M$ are involutive.
\end{lem}
\begin{proof}
For $e_1, e_2 \in \Gamma(\mathcal{E}M)$ and $f_1, f_2 \in \Gamma(\mathcal{F}M)$, from $\nabla P =0$, we get
\begin{align}
\label{EtimesEinE} P(\nabla_{e_1} e_2) =& \nabla_{e_1} P(e_2) = \nabla_{e_1} e_2,\\
\label{FtimesFinF} P(\nabla_{f_1} f_2) =& \nabla_{f_1} P(f_2) = \nabla_{f_1} 0 =0.
\end{align}
i.e. $\nabla_{e_1} e_2 \in \Gamma(\mathcal{E} M) $ and $\nabla_{f_1} f_2 \in \Gamma(\mathcal{F} M)$. As the Levi-Civita connection is torsion free we have
that $\mathcal{E} M$ and $\mathcal{F} M$ are involutive.
\end{proof}

By the Theorem of Frobenius
, these distributions define two foliations, such that, at each $p \in M$, there are two leaves $\mathcal{E}_p$ and $\mathcal{F}_p$ that intersect orthogonally at $p$. We want to understand what they are. The inclusions $i_{\mathcal{E}_p}$ and $i_{\mathcal{F}_p}$ of these leaves are immersions:
\begin{table}[ht]
\begin{minipage}[b]{0.45\linewidth}\centering
\begin{displaymath}
   \xymatrix{M \ar[r]^{F} & \R^n \\
        \mathcal{E}_p  \ar[u]^{i_{\mathcal{E}_p}} \ar[ur]_{F \circ i_{\mathcal{E}_p}} & }
\end{displaymath}
\end{minipage}
\begin{minipage}[b]{0.45\linewidth}
\centering
\begin{displaymath}
    \xymatrix{M \ar[r]^{F} & \R^n \\
        \mathcal{F}_p \ar[u]^{i_{\mathcal{F}_p}} \ar[ur]_{F \circ i_{\mathcal{F}_p}} & }
\end{displaymath}
\end{minipage}
\end{table}

We need the symmetric (by Lemma \ref{equationswithindices}) tensor
$$
(P*A)_{ij}:=P^k_iA_{kj}.
$$

\begin{lem}\label{mugujininho}
Under the hypothesis of Theorem \ref{ResultNoncompact1}, the following equations hold:
\begin{align}
\label{ThetaNablaA} \theta^k \nabla^\perp_k A_{ij} = 0,\\
\label{A,PA} A_{ij}=P_i^kA_{kj}.
\end{align}
\end{lem}

\begin{proof}
First, from \eqref{secondderivativeH} (with $\nabla \H =0$) and \eqref{P,P2}, we get
\begin{equation}
\label{ThetaNablaPAequalsZero} \theta^k\nabla^\perp_k (P_i^lA_{lj}) = P_i^lA_{lj} - P_i^lA_{lj} = 0.
\end{equation}
To prove \eqref{A,PA}, it is enough to show 
$$
\|A_\pm\|^2 = \|P*A_\pm\|^2.
$$
One sees this using eq. \eqref{P,P2} to calculate $\|A_\pm - P*A_\pm\|^2 = \|A_\pm\|^2 - \|P*A_\pm\|^2$.

Let us then prove that $\|A_\pm\|^2 = \|P*A_\pm\|^2$. First of all, using eq. \eqref{secondderivativeH} (with $\nabla_i \H =0$), eq. \eqref{P,P2} and $\theta^k\nabla^\perp_k (P_i^lA_{lj}) =0$, it holds that
\begin{equation}
\label{AequalPA}\theta^k\nabla_k(\|A_\pm\|^2 - \|P*A_\pm\|^2)
= -2 (\|A_\pm\|^2 - \|P*A_\pm\|^2).
\end{equation}
If $\theta =0$ at some point $p \in M$, then this equation implies $\|A_\pm\|^2 =\|P*A_\pm\|^2$ and $A_{ij}=P_i^l A_{lj}$ at this point. So, without loss of generality, we can consider only the points $q \in M$ with $\theta (q) \neq 0$. Fix one of these and consider the integral curve $\gamma:(-a,b) \to M$ of $\theta$
with $\gamma(0)= q$, for some $a,b > 0$. Along this curve we define the function
$$
f(s):= \|\theta\|^2(\gamma(s)),
$$
and get
$$
\frac{d}{ds}f = \nabla_{\dot{\gamma}} \|\theta\|^2 = \theta^k \nabla_k \|\theta\|^2 = 2\theta^k \theta^l \nabla_k \theta_l.
$$
but, from $\H =-F^\perp$,
$$
\nabla_i \theta_j = \nabla_i \langle F, F_j \rangle = g_{ij} - \langle \H, A_{ij} \rangle
$$
and $\theta^i \nabla_i \theta_j = \theta_j$ because of equation \eqref{thetaA}, so that
\begin{equation}\label{Integralcurve1}
\frac{d}{ds}f=2f.
\end{equation}
This has a unique solution with $f(0)=\|\theta\|^2(q)$:
$$
\|\theta\|^2 (\gamma(s)) = \|\theta\|^2(q) e^{2s}>0,
$$
in particular $\|\theta\|^2 (\gamma(s)) \neq 0$ for all $s \in (a,b)$, then these integral curves do not cross any singular point and the maximal integral curve is defined for all $\R$  and it is not closed (because of injectivity of $e^{2s}$). Over this same curve we define functions $\tilde{f_\pm}: \R \to \R$,
$$
\tilde{f_\pm}(s):= (\|A_\pm\|^2 - \|P*A_\pm\|^2)(\gamma(s))
$$
and, using equation \eqref{AequalPA}, get $\frac{d\tilde{f_\pm}}{ds} = -2\tilde{f_\pm}$.
This has a unique solution with $\tilde{f_\pm}(0) = (\|A_\pm\|^2 - \|P*A_\pm\|^2)(q)$:
$$
(\|A_\pm\|^2 - \|P*A_\pm\|^2)(\gamma(s)) = (\|A_\pm\|^2 - \|P*A_\pm\|^2)(q)e^{-2t}.
$$
If $(\|A_\pm\|^2 - \|P*A_\pm\|^2)(q) \neq 0$, then $(\|A_\pm\|^2 - \|P*A_\pm\|^2)(\gamma(s)) \to \pm \infty$ as $s \to -\infty$ and this contradicts the boundedness of $\|A_\pm\|^2$. So $A= P*A$ and \eqref{ThetaNablaPAequalsZero} implies \eqref{ThetaNablaA}.
\end{proof}
Let us now examine the leaves of the distribution $\mathcal{E}M$.
\begin{lem}\label{jnjnjnjnkmkmkmkmk}
Under the hypothesis of Theorem \ref{ResultNoncompact1} it holds that $\mathcal{E}_p$ is immersed into $\mathcal{H}^{n-1}(\|F\|^2(p))$ through $F \circ i_{\mathcal{E}_p}$. Beyond this, $\mathcal{E}_p$ is geodesically complete and there is $q \in M$ so that $F \circ i_{\mathcal{E}_q}$ is a minimal immersion into $\mathcal{H}^{n-1}(\|F\|^2(q))$.
\end{lem}
\begin{proof}
$\mathcal{E}_p$ is an $r$-dimensional manifold immersed in $M$ under the natural inclusion $i_{\mathcal{E}}$. Let $A_{F \circ i_\mathcal{E}}$ and $A_{i_\mathcal{E}}$ denote the second fundamental tensors of $F \circ i_{\mathcal{E}}$ and $i_{\mathcal{E}}$.

From equation \eqref{chainrulesecondfund} 
it holds that
\begin{equation}\label{AfiequalAfplusDfAi}
A_{F \circ i_\mathcal{E}} = A_F + dF (A_{i_\mathcal{E}}).
\end{equation}
On the other hand one can write, for local vector fields $e_1, e_2 \in \Gamma(T \mathcal{E}_p)$,
\begin{equation} \label{seconddefinitionAE}
A_{i_\mathcal{E}}(e_1,e_2) = \nabla_{e_1} e_2 - \nabla'_{e_1} e_2,
\end{equation}
where $\nabla'$ is the Levi-Civita connection of $\mathcal{E}_p$ (with respect to the induced metric). But $\nabla_{e_1} e_2 \in \Gamma(\mathcal{E} M)$ by eq. \eqref{EtimesEinE} and $di_\mathcal{E}(\nabla'_{e_1} e_2) \in \Gamma(\mathcal{E}M)$, yet $A_{i_\mathcal{E}}(e_1,e_2) \in \Gamma(T\mathcal{E}^\perp)$, so
\begin{equation}\label{popopopopopopopopopopo}
A_{i_\mathcal{E}}=0.
\end{equation}
Then, in particular, the geodesics of $\mathcal{E}_p$ are also geodesics of $M$ and, 
as $M$ is geodesically complete, so is $\mathcal{E}_p$.

From equation \eqref{kernelTheta} we get, for any $q\in \mathcal{E}_p$ and all $V \in \mathcal{E}M_{i_\mathcal{E}(q)}$, $V=V^i\frac{\partial}{\partial x^i}$,
$$
0=\theta_j V^j = \langle F, F_j\rangle V^j = \langle F, dF(V)\rangle,
$$ 
which means that $F (i_\mathcal{E}(q))\in T_q\mathcal{E}_p^\perp$ (the normal space of $F \circ i_\mathcal{E}$ at $q$) and
\begin{equation}\label{srewsrewrewrewerw}
V \|F\|^2 = 2V^j \langle F_j, F \rangle = 2\langle dF(V), F \rangle =0, 
\end{equation}
so that $\|F\|^2$ is constant on the leaf $\mathcal{E}_p$ (but it depends on $p$), and $\mathcal{E}_p$ is immersed, through $F \circ i_\mathcal{E}$, in the hyperquadric $\mathcal{H}^{n-1}(\|F\|^2(p))$. 

Let us now take a look at a special leaf of the distribution $\mathcal{E}M$. Because of Remark \ref{seilaoque} 
and Lemma \ref{LemmaDoido} 
there is a point $q \in M$, with $\|F(q)\|^2 = \min_{x \in M} {\|F(x)\|^2}$. Let us consider the leaf $\mathcal{E}_q$. We are showing that $F(M)$ is some cylinder and figure \ref{mainlypositive11} shows the intersection of a cylinder with two spheres. The small circle in the middle is a minimal surface of the smallest sphere but the two other circles are not a minimal surfaces of the bigger sphere.
\begin{figure}[h]
\centering\includegraphics[height=53mm, angle=-90]{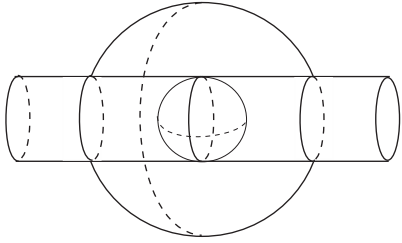}
\caption{Intersection of a cylinder with spheres}\label{mainlypositive11}
\end{figure}

The norm of $F$ must be constant over this leaf by eq. \eqref{srewsrewrewrewerw}, so that all the points of the leaf minimize the norm of $F$. But then, $2 \langle dF(X), F \rangle=X \|F\|^2 = 0$ for any $X \in T_{q'}M$, $q' \in i_\mathcal{E}(\mathcal{E}_q)$, this means that $F(q')$ is orthogonal to $T_{q'}M$, i. e.
\begin{equation}\label{FperpequalsFfirstcase}
F^\perp(q') =F(q').
\end{equation}

We claim that $F \circ i_\mathcal{E} (\mathcal{E}_q)$ is a minimal surface of the hyperquadric $\mathcal{H}^{n-1}(\|F\|^2(q))$. First, the Levi-Civita connection on the hyperquadric is given by the projection $(\Pr_{\mathcal{H}^{n-1}})$ of the Levi-Civita connection of $\R^{q,n}$, which we denote $D$, over the tangent bundle of the hyperquadric. Then, using eq. \eqref{AfiequalAfplusDfAi} with $A_{i_\mathcal{E}} = 0$, it holds 
\begin{equation}
\label{lkmlkmlkmlkmlksffasdsf} A_{\mathcal{H}^{n-1}} (X,Y)=\Pr_{\mathcal{H}^{n-1}} (D_X Y) - \nabla'_X Y = \Pr_{\mathcal{H}^{n-1}} (A_{F \circ i_\mathcal{E}}(X,Y)) = \Pr_{\mathcal{H}^{n-1}} (A_F(X,Y))
\end{equation}

On the other hand, take a vector $V \in T_{q'}M^\perp$, $q' \in \mathcal{E}_q$, 
then, using that $P^{ij}A_{ij}$ is in the same direction as $\H$ (eq. \eqref{PAdirectionNU}), one gets $\langle P^ {ij} A_{ij}, V \rangle = 0.$

Take an orthonormal basis $\{e_1, \ldots, e_r, f_1, \ldots, f_{m-r}\}$ of $T_{q'}M$ such that $\{e_1, \ldots, e_r\}$ is a basis of $\mathcal{E} M_{q'}$ and $\{f_1, \ldots, f_{m-r}\}$ is a basis of $\mathcal{F} M_{q'}$, then
$$
\mbox{tr}_\mathcal{E} \langle A, V\rangle = \sum_{i=1}^r \langle V, A(P(e_i), e_i) + \sum_{i=1}^{m-r} \langle V, A(P(f_i), f_i) \rangle = \mbox{tr}_M \langle V, P*A \rangle = 0,
$$
where we used that $P(e_i) = e_i$ and $P(f_i) = 0$. This holds for any $q' \in i_\mathcal{E}(\mathcal{E}_q)$ and means that $\mbox{tr}_\mathcal{E} A = \mathfrak{a}(x) \H$ for some continuous function $\mathfrak{a}:\mathcal{E}_q \to \R$. By eq. \eqref{FperpequalsFfirstcase}, eq. \eqref{lkmlkmlkmlkmlksffasdsf} and denoting $\H_{\mathcal{H}^{n-1}}$ the mean curvature vector of the immersion of $\mathcal{E}_q$ into $\mathcal{H}^{n-1}(\|F\|^2(q))$, we get at $q'$
$$
\H_{\mathcal{H}^{n-1}} = \Pr_{\mathcal{H}^{n-1}} (\mbox{tr}_\mathcal{E} A) = \Pr_{\mathcal{H}^{n-1}} (\mathfrak{a} \H) = \Pr_{\mathcal{H}^{n-1}} (-\mathfrak{a} F^\perp) = \Pr_{\mathcal{H}^{n-1}} (-\mathfrak{a}(x) F) = 0,
$$
because the position vector is orthogonal to the hyperquadric. Then $\mathcal{E}_q$ is a minimal surface of the hyperquadric $\mathcal{H}^{n-1}(r)$, because $\|F\|^2(q) = \|\H \|^2(q) =r$ by eq. \eqref{normHsquaredequalr}.
\end{proof}

We will now analyse the leaves of the distribution $\mathcal{F}M$.

\begin{lem}\label{jujujujujujjujujujujujjujuj}
Under the hypothesis of Theorem \ref{ResultNoncompact1}, it holds that $F\circ i_\F(\mathcal{F}_p)$ is an affine space in $\R^{q,n}$ for any $p \in M$. Beyond that, if $q \in i_\mathcal{E}(\mathcal{E}_p)$, then $F\circ i_\F(\F_p)$ and $F\circ i_\F(\F_q)$ are parallel.
\end{lem}
\begin{proof}
First, we show that they are affine subspaces of $\R^{q,n}$. Let $q \in i_\F (\mathcal{F}_p)$ be an arbitrary point and $f \in \mathcal{F}M_q$ and $X \in T_qM$ be vectors, then eq. \eqref{A,PA} implies
\begin{equation} \label{pocoyo}
A(f,X) = X^j f^i A_{ij} = X^j f^i P_i^kA_{kj} = 0
\end{equation}
because $\F M_q$ is formed by the vectors in the null space of $P$.

Let us denote $A_{F \circ i_\F}$ and $A_{i_\F}$ the second fundamental tensors of the immersions $F \circ i_{\F}$ and $i_{\F}$ respectively. From equation \eqref{chainrulesecondfund}:
$$
A_{F \circ i_\F} = A_F + dF (A_{i_\F}).
$$
Equation \eqref{pocoyo} implies that $A_{F \circ i_\F}=dF (A_{i_\F})$. Thence 
\begin{equation}
A_{F \circ i_\F}=A_{i_\F} = 0,
\end{equation}
which means that $i_\F$ is totally geodesic. 
One has $D_{f_1} f_2 = A_F(f_1, f_2) + \nabla_{f_1} f_2 = \nabla_{f_1} f_2$ so that the geodesics of $\F_p$ are also geodesics of $\R^{q,n}$ (which are straight lines). Furthermore  $\F_p$ is geodesically complete, 
so that each connected component of $F \circ i_\F(\F_p)$ is an affine $m-r$-dimensional subspace of $\R^{q,n}$ by the linearity of $d(F \circ i_\F)$, 
but $F \circ i_\F(\F_p)$ is a leaf and thence a whole affine $m-r$-dimensional subspace of $\R^{q,n}$.

Let us fix $p$ and prove that $F \circ i_\F(\F_{p'})$ is parallel to $F \circ i_\F(\F_p)$ for any $p' \in i_\mathcal{E}(\mathcal{E}_p)$. For this, let $\gamma: [0,1]\to i_\mathcal{E}(\mathcal{E}_p)$ be a smooth curve with $\gamma(0) =p$ and $\gamma(1)=p'$ and let $f_p \in \F M_p$ be a vector. Denote by $f(t)$ the parallel transport of $f_p$ along $\gamma$ with respect to $\nabla$. From $\nabla P = 0$, it holds that $\nabla_{\dot{\gamma}} (P f(t))= P(\nabla_{\dot{\gamma}} f(t))=0$, which means that $P (f(t))$ is the parallel translation of $P(f_p)=0$, so that $P (f(t))=0$ and $f(t) \in \F M_{i_\mathcal{E}\circ \gamma(t)}$ for all $t \in [0,1]$. By using eq. \eqref{pocoyo},
$$
D_{\dot{\gamma}}f(t) = \nabla_{\dot{\gamma}} f(t) + A(\dot{\gamma}, f(t)) = 0,
$$
where $D$ is the Levi-Civita connection of $\R^{q,n}$. This means that $dF \circ f(t)$ is also the parallel translation of $dF \circ f_p$ in $\R^{q,n}$, which is $dF(f_p)$ for all $t \in [0,1]$. This means that $dF \circ di_\F (T_p\F_p) \subset dF \circ di_\F(T_{p'}\F_p)$. Analogously, $dF \circ di_\F(T_{p'}\F_p) \subset dF \circ di_\F(T_p\F_p)$. But we already know that the leaves of $\F$ are affine subspaces, so that they are equal, up to a translation, to their tangent spaces and thus are parallel.
\end{proof}

All that is left of Theorem \ref{ResultNoncompact1} is to show that $F(M)$ is the product $F(\mathcal{E}_q) \times F(\F_q)$, where $q \in M$ minimizes $\|F\|^2$.

Let $q \in M$ be a minimal point of $\|F\|^2$ and $\{f_1, \ldots, f_{m-r}\}$ be an orthonormal basis of $\F M_q$. We define a function $\mathfrak{h}: \mathcal{E}_q \times \R^{m-r} \to F(M)$, given by
$$
\mathfrak{h}(p,X) = F(i_{\mathcal{E}}(p)) + X^i dF(f_i) \hspace{20mm} \, \forall X= (X^1, \ldots, X^{m-r}) \in \R^{m-r}, \, p\in \mathcal{E}_q.
$$
As all the leaves $\F_{q'}$, $q' \in \mathcal{E}_q$, are parallel, the image of $\mathfrak{h}$ is indeed contained in $F(M)$. Let us consider in $\R^{m-r}$ the canonical metric and in $\mathcal{E}_q \times \R^{m-r}$ the product metric, so that $\mathfrak{h}$ is an isometry because $F$ and $i_\mathcal{E}$ are isometries.

$\mathcal{E}_q \times \R^{m-r}$ is geodesically complete. We claim that $\mathfrak{h}$ is surjective. To see this, take $(p,X) \in \mathcal{E}_q \times \R^{m-r}$, $y:=\mathfrak{h}(p,X)\in F(M)$, and $z \in F(M)$. Let $y'\in M$ and $z' \in M$ be such that $F(y')=y$ and $F(z')=z$. From the fact that $M$ is geodesically complete, there is a vector $Y \in T_{y'} M$ such that $\exp (Y) =z'$ (by the Theorem of Hopf and Rinow). Then decompose $Y =Y_1+ Y_2$ with $Y_1 \in T_p\mathcal{E}_q$ and $Y_2 = Y_2^l f_l(p) \in \mathcal{F}M_p$. Now denote $Y_{20}:=(Y_2^1, \ldots, Y_2^{m-r})$, then for the exponential in $\mathcal{E}_q \times \R^{m-r}$ it holds that
\begin{align*}
\mathfrak{h} (\exp (Y_1, Y_{20})) =& \exp (d\mathfrak{h}(Y_1, Y_{20})) = \exp (dF \circ di_\mathcal{E} (Y_1) + dF (Y_2)) \\
=& F (\exp( di_{\mathcal{E}}(Y_1) +Y_2)) = F(\exp (Y))=z,
\end{align*}
where we understand $F(M)$ locally as a manifold (isometric to $M$ and with the same dimension) and thence define the exponential there locally, so that, by the compactness of the domain of the geodesic segment connecting $y'$ and $z'$, the exponential is well defined. This proves that $z \in \mathfrak{h}(\mathcal{E}_q \times \R^{m-r})$.

Then $F(M)$ is the product of an affine space with a minimal surface of the hyperquadric $\mathcal{H}^{n-1}(r)$ with $\|\H \|^2 =r$.
\end{proof}
\begin{rem}
The induced (from $\R^{q,n}$) inner product on the affine space has to be positive definite, because $F$ is spacelike.
\end{rem}

\end{section}

\begin{section}{The Second Case}
%
%
\begin{thm}\label{ResultNoncompact2}
Let $M$ be a smooth manifold and $F:M \to \R^{q,n}$ be a mainly positive, spacelike, shrinking self-similar solution of the mean curvature flow with bounded geometry such that $F(M)$ is unbounded. Beyond that, let $F$ satisfy the conditions:  $\|\H\|^2(p) \neq 0$ for all $p \in M$ and the principal normal is parallel in the normal bundle ($\nabla^\perp \nu \equiv 0$). If $\nabla \|\H\|(p) \neq 0$ for some $p \in M$, then
\begin{equation}\label{result*}
F(M) = \Gamma \times \R^{m-1},
\end{equation}
where $\Gamma$ is a rescaling of an Abresch \& Langer curve in a spacelike plane and $R^{m-1}$ is an $m-1$ dimensional spacelike affine space in $\R^{q,n}$.
\end{thm}
\begin{proof}
Let $p \in M$ be a point with $\nabla \| \H \|(p) \neq 0$ and $\{e_i\}_{1=1,\ldots,m}$ an orthonormal basis of $T_pM$ made by the eigenvectors of $P$ with $e_1 =\frac{\nabla \|\H\|}{\|\nabla \|\H\|\|}(p)$. From Lemma \ref{cucuxico} $P$ (in this basis) has only one nonzero element. But $\mbox{tr}(P) = \sum_i \langle \H, A(e_i, e_i) \rangle = \|\H\|^2$, so $\|\H\|^2$ is the eigenvalue associated with $\frac{\nabla \|\H\|}{\|\nabla \|\H\|\|}$:
\begin{equation} \label{EigenvalueP}
P^i_j\frac{\nabla^j \|\H\|}{\|\nabla \|\H\|\|} = \|\H\|^2 \frac{\nabla^i \|\H\|}{\|\nabla \|\H\|\|}
\end{equation}
and $\| P \|^2 = \| \H \|^4$ at this point, but this equation holds all over $M$ because equation \eqref{importantintegral} together with Lemma \ref{lemaImportante} implies
$$
\nabla \left(\frac{\|P\|^2}{\|\H\|^4} \right)=0.
$$

\begin{rem}\label{queroirembora}
Let us choose Riemannian normal coordinates on a neighborhood of $p$ such that $\frac{\partial}{\partial x^i}(p) = e_i$, then it holds in $p$: $P_{ij} = \|\H\|^2 \delta_{1i}\delta_{1j}$ and $g_{ij}= \delta_{ij}$, thus it follows that $P_i^k A_{kj} = 0$ if $i \neq 1$ and, from Lemma \ref{equationswithindices}, item (3), that $P_i^k A_{kj} = 0$ if  $j \neq 1$ so that $P_1^k A_k^1 = P_i^k A_k^i$, which is in the direction of $\nu/\|\H\|$ from Lemma \ref{equationswithindices} (item 1), thence the component of $P_i^k A_{kj}$ in $F^*(\R^{q,n})$ has the same direction as $\frac{\nu}{\|\H\|}$ for $i,j \in 1, \ldots, m$ and
\begin{equation}\label{PAPnu}
P_i^k A_{kj} = P_i^k \langle \nu, A_{kj} \rangle \nu = \frac{1}{\|\H\|} P_i^k P_{kj} \nu = \|\H\|^3 \delta_{1i}\delta_1^k \delta_{1k}\delta_{1j}\nu = \|\H\| P_{ij}\nu,
\end{equation}
i. e. $P*A = \|\H\| P \otimes \nu $, but this can be done for any $p \in M$ with $\nabla \|\H\|(p) \neq 0$, so that this tensor equality holds in any region of $M$ that satisfies $\nabla \|\H\|(p) \neq 0$.
\end{rem}

Let us now define
$$
\mathring{M} = \{p \in M: \nabla \|\H\| \neq 0\}.
$$
Which is 
open, and thus a submanifold of $M$ (possibly incomplete). Let $U \subset M$ be a connected component of $\mathring{M}$. We take, over $U$, the distributions $\mathcal{E}U$ and $\mathcal{F}U$ given by
\begin{align}
\mathcal{E}U_p:=&\{V \in T_p U: PV= \|\H\|^2 V\},\\
\label{kernFU}\mathcal{F}U_p:&=\{V \in T_p U: PV= 0\}.
\end{align}



In order to investigate these distributions we need further information about the second funcamental tensor. For this purpose we define the tensor $\mathring{A}:= A - \frac{1}{\|\H\|}P \otimes \nu$.

\begin{lem}\label{voumeembora}
Under the hypothesis of Theorem \ref{ResultNoncompact2}, it holds
\vspace{-1mm}
\begin{equation}
 \label{ANormHPnu} A_{ij}= \frac{1}{\|\H\|}P_{ij}\nu.
\end{equation}
\end{lem}
\begin{proof}
From eq. \eqref{HThetaA} we calculate
\begin{equation}\label{HnablaHthetaP}
\|\H\| \nabla_i \|\H\| = \langle \nabla_i^\perp \H, \H \rangle = \theta^k P_{ik},
\end{equation}
this equation together with eq. \eqref{EigenvalueP} 
implies
\begin{equation}\label{ThetaNablaHNablaH2}
\theta (\nabla \|\H\|) = \frac{\theta^k}{\|\H\|^2} P_k^l (\nabla_l \|\H\|) = \frac{\|\H\|\nabla^l \|\H\|\nabla_l \|\H\|}{\|\H\|^2} = \frac{\|\nabla\|\H\|\|^2}{\|\H\|}.
\end{equation}
In order to attain $\mathring{A}_{ij} = 0$ we consider the integral curves of the projection of $F$ in $\F U$:
$$
\mathring{\theta}_i = \theta_i - \frac{\theta(\nabla^k\|\H\|)}{\|\nabla_k\|\H\|\|^2}\nabla_i\|\H\| = \theta_i - \frac{1}{\|\H\|}\nabla_i\|\H\|.
$$
For $\mathring{A}$ one calculates
$$
\nonumber\mathring{\theta}^k\mathring{A}_{ki} = \theta^k A_{ki} - \nabla_i^\perp \H = 0
$$
using eqs. \eqref{EigenvalueP}, \eqref{PAPnu}, \eqref{HnablaHthetaP} and \eqref{HThetaA}. From this
\begin{align}
\nonumber 0 = \nabla_l^\perp (\mathring{\theta}^k \mathring{A}_{ki}) =& \left(\delta_l^k - P_l^k + \frac{1}{\|\H\|^2}\nabla_l\|\H\|\nabla^k\|\H\| - \frac{1}{\|\H\|}\nabla_l \nabla^k \|\H\| \right)\mathring{A}_{ki} + \mathring{\theta}^k\nabla^\perp_l \mathring{A}_{ki}\\
\label{firststep} 0 =&\mathring{A}_{li} -\frac{1}{\|\H\|}\nabla_l \nabla^k \|\H\|\mathring{A}_{ki} + \mathring{\theta}^k\nabla^\perp_l \mathring{A}_{ki},
\end{align}
because equations \eqref{EigenvalueP} and \eqref{PAPnu} with $P^k_l P_{ki} = \|\H\|^2 P_{li}$ imply that $P_l^k \mathring{A}_{ki} = 0$ and equations \eqref{EigenvalueP}, \eqref{PAPnu} and \eqref{HnablaHthetaP} imply
\begin{align*}
\nabla^k\|\H\| \mathring{A}_{ki} =& \nabla^k\|\H\| A_{ki} - \nabla^k \|\H\| \frac{1}{\|\H\|}P_{ki} \nu =\frac{1}{\|\H\|}\theta^l P_l^k A_{ki}-\|\H\| \nabla_i\|\H\| \nu\\
=&\theta^l P_{li} \nu -\|\H\| \nabla_i\|\H\| \nu = \|\H\| \nabla_i\|\H\| \nu -\|\H\| \nabla_i\|\H\| \nu=0.
\end{align*}
\indent On the other hand, using the Riemannian normal coordinates of Remark \ref{queroirembora} we calculate for $i \in 1,\ldots,n$ and $j \in 2, \ldots, n$,
$$
\nabla_i \nabla_j \|\H\| = \frac{\partial}{\partial x^i} \left< \nabla \|\H\|, \frac{\partial}{\partial x^j} \right> - \left< \nabla\|\H\|, \nabla_{\frac{\partial}{\partial x^i}} \frac{\partial}{\partial x^j} \right> = 0
$$
because $\frac{\partial}{\partial x^j} \in \F U_p$ and $\nabla_{\frac{\partial}{\partial x^i}} \frac{\partial}{\partial x^j} (p) = 0$. But $\nabla_i \nabla_j \|\H\|$ is symmetric. So that $\nabla_i \nabla_j \|\H\|$ is nonzero only if $i=j=1$, so that using $P_{ij}= \|\H\|^2 \delta_{1i}\delta_{1j}$ we get
\begin{equation}\label{LaplaceNablaNabla}
\nabla_i \nabla_j \|\H\| = \frac{\triangle \|\H\|}{\|\H\|^2} P_{ij},
\end{equation}
but this equation is tensorial 
then holds for any coordinate map. And, as we could do the same for every point $p \in U$, this holds in the whole $\mathring{M}$. So that, with eq. \eqref{PAPnu},
\begin{align*}
\nabla_l \nabla^k \|\H\|\mathring{A}_{ki} =& \frac{\bigtriangleup\|\H\|}{\|H\|^2}P_l^k\mathring{A}_{ki}
= \frac{\bigtriangleup\|\H\|}{\|H\|^2}\left( \|\H\| P_{li}\nu - \|\H\|P_{li}\nu \right) = 0
\end{align*}
and equation \eqref{firststep} turns out to be
\begin{equation}\label{TethaNablaAA}
\mathring{\theta}^k\nabla^\perp_l \mathring{A}_{ki} = -\mathring{A}_{li}.
\end{equation}
\vspace{-2mm}Finally, we calculate using eq. \eqref{HnablaHthetaP}
\begin{align*}
\nabla_i\|\mathring{\theta}\|^2 = 2\nabla_i \mathring{\theta}_l \mathring{\theta}^l=& 2\left(\theta_i -\theta^l P_{li} + \theta^l\frac{\nabla_i\|\H\| \nabla_l\|\H\|}{\|\H\|^2} - \frac{\triangle\|\H\|}{\|\H\|^3}P_{il}\theta^l -\frac{\nabla_i\|\H\|}{\|\H\|}\right.\\
&+ \left.\frac{\nabla^l\|\H\|}{\|\H\|}P_{li} - \frac{\|\nabla\|\H\|\|^2}{\|\H\|^3}\nabla_i\|\H\| + \frac{\triangle \|\H\|}{\|\H\|^4}P_{il}\nabla^l\|\H\| \right)\\
=& 2\left(\theta_i + \theta^l\frac{\nabla_i\|\H\| \nabla_l\|\H\|}{\|\H\|^2} -\frac{\nabla_i\|\H\|}{\|\H\|}- \frac{\|\nabla\|\H\|\|^2}{\|\H\|^3}\nabla_i\|\H\| \right),
\end{align*}
but using eq. \eqref{HnablaHthetaP}
\begin{align*}
\mathring{\theta}^i \nabla_i \|\H\| &
= \langle F, F^i \rangle \nabla_i\|\H\| -\langle F, F^k \rangle \nabla_k \|\H\|=0.
\end{align*}
which was already expected
, so that it holds

\begin{equation}\label{ThetaNablaNormTheta}
\mathring{\theta}(\nabla \|\mathring{\theta}\|^2) = 2\mathring{\theta}^i \theta_i = 2\mathring{\theta}^i \left(\theta_i -\frac{\nabla_i\|\H\|}{\|\H\|} \right) = 2\|\mathring{\theta}\|^2.
\end{equation}

We follow then as in Lemma \ref{mugujininho} but still have to check if eq. \eqref{ThetaNablaNormTheta} holds in the whole $M$. In open sets of $M \setminus \mathring{M}$, the equations of the first case hold. $P=P^2$ (eq. \eqref{P,P2}) implies that the only non-zero eigenvalue of $P$ is 1, then $\nabla P = 0$ together with $\|P\|^2=r$, where $r$ is the multiplicity of the eigenvalue $1$, implies that $r$ is constant, but $M$ is connected and, in $\overline{\mathring{M}}$, the tensor $P$ has only one non-zero eigenvalue and it has multiplicity 1, then $P$ has only one non-zero eigenvalue and it has multiplicity 1 also in open sets of $M \setminus \mathring{M}$. Therefore, as in Remark \ref{queroirembora}, $P_{ik} A^k_j$ is in the direction of $\frac{\H}{\|\H\|}$ and
\begin{align}
\nonumber \|\mathring{A}_\pm\|^2 =& \left\|A_\pm - \frac{1}{\|\H\|}P \otimes \nu_\pm\right\|^2= \|A_\pm\|^2 -\frac{2}{\|H\|}\langle A^{ij}_\pm, P_{ij}\nu_\pm \rangle +\frac{1}{\|\H\|^4}\|P\|^2\|\H_\pm\|^2\\
\nonumber =&\|A_\pm\|^2 -\frac{2}{\|H\|^4}\langle P^{ij}P_{ij}\H_\pm, \H_\pm \rangle +\frac{1}{\|\H\|^4}\|P\|^2\|\H_\pm\|^2\\
\label{indoembora} \|\mathring{A}_\pm\|^2=&\|A_\pm\|^2 -\frac{1}{\|\H\|^4}\|P\|^2\|\H_\pm\|^2,
\end{align}
and
\begin{align*}
\| P*A_\pm\|^2 =& \langle P_{ik} (A_\pm)^k_j, P_l^i A_\pm^{lj} \rangle = \left< P_{ik} \left< A^k_j, \frac{\H}{\|\H\|} \right> \frac{\H_\pm}{\|\H\|},P_l^i \left< A^{lj}, \frac{\H}{\|\H\|} \right> \frac{\H_\pm}{\|\H\|}\right>\\
=& P_{ik} P^k_j P_l^i P^{lj}\frac{1}{\|\H\|^4}\|\H_\pm\|^2 = \frac{\|P\|^2}{\|\H\|^4}\|\H_\pm\|^2,
\end{align*}
using $P_{ij} = P_{ik} P^k_j$ (eq. \eqref{P,P2}), thence, by eq. \eqref{indoembora}, it holds
\begin{equation}
\label{AringAPA}\|\mathring{A}_\pm\|^2= \|A_\pm\|^2 -\frac{1}{\|\H\|^4}\|P\|^2\|\H_\pm\|^2 = \|A_\pm\|^2 - \| P*A_\pm\|^2
\end{equation}
and equation $\eqref{AequalPA}$ implies that $\frac{d\tilde{f}}{ds} = -2\tilde{f}$ holds in the whole manifold $M$. This O.D.E. has a unique solution $\|\mathring{A}\|^2(\gamma(s)) = \|\mathring{A}\|^2(q)e^{-2s}$. If $\|\mathring{A}_\pm\|^2(q) \neq 0$ then $\|\mathring{A}_\pm\|^2(\gamma(s)) \to \pm\infty$ as $s \to -\infty$ and this contradicts the boundedness of the second fundamental tensor (by the definition of bounded geometry). So that $\|\mathring{A}_\pm\|^2=0$ as we wanted to show.
\end{proof}

\begin{lem}
Under the hypothesis of Theorem \ref{ResultNoncompact2} the distributions $\mathcal{E}U$ and $\mathcal{F}U$ are involutive.
\end{lem}

\begin{proof}
First of all recall that $\mathcal{E} U$ is spanned by $\overrightarrow{e}:=\nabla \|\H\| /\|\nabla \|\H\|\|(p)$ at any $p \in U$. 
%
Then let $p\in U$ be a point and $V_p \in T_p U$ be normal to $\nabla \|\H\| /\|\nabla \|\H\|\|(p)$; beyond this, let $V \in \Gamma(\left. T U\right|_\Omega)$ be the parallel transport of $V_p$ over all geodesics through $p$ in a small neighborhood $\Omega$ of $p$. 
Then, as any $X,Y \in \Gamma(\mathcal{E}U)$ are of the form $X = x\overrightarrow{e}$ and $Y = y\overrightarrow{e}$ for some $x, y \in C^1 (U)$, it holds at $p$ that
\begin{equation}\label{PnablaXYinEU}
\left< \nabla_{x\overrightarrow{e}} y\overrightarrow{e}, V \right> = x\overrightarrow{e}\left<y\overrightarrow{e}, V \right> - \left< y\overrightarrow{e}, \nabla_{x\overrightarrow{e}} V \right>=0.
\end{equation}
But this could be done for any $p \in U$, so that, just as in the first case, $\mathcal{E}U$ is in particular involutive. So, by the Theorem of Frobenius, there is a foliation, whose tangent spaces of the leaves are given by this distribution. The leaves are immersed in $M$ and, again as in the first case, they are totally geodesic (analogous to eq. \eqref{seconddefinitionAE} and eq. \eqref{popopopopopopopopopopo}). This means, in particular, for any $p \in U$, that a geodesic of the one dimensional leaf ($\mathcal{E}_p$) that goes through $p$ is also a geodesic of $U$.

%
%
Let $p \in U$ be a fixed point and take Riemannian normal coordinates around $p$ such that $\frac{\partial}{\partial x^1} = \frac{\nabla \|\H\|}{\|\nabla \|\H\|\|}$. This way the tensor $P$ is written, in these coordinates, as $P_{ij} = \|\H\|^2 \delta_{1i}\delta_{1j}$. So, for $V,W \in \Gamma(\F U)$, using the fact that $\F U_p \perp \mathcal{E} U_p$, we get
$$
0=\langle \nabla \|\H\|, V \rangle = \nabla_V \|\H\|,
$$
from this, remembering $\Gamma_{ij}^k(p)=0$, follows
$$
\nabla_V P = \nabla_V \|\H\|^2 \delta_{1i}\delta_{1j}= 0 \,\,\, \mbox{and} \,\,\, P(\nabla_V W)= \nabla_V (PW) = 0.
$$

This means
\begin{equation}\label{FtimesFinF2}
\nabla_V W \in \Gamma(\F U) \,\,\forall\,\, V,W \in \Gamma(\F U).
\end{equation}
As this holds for any $p \in U$ and the final expressions do not depend on local coordinates this holds in the whole $\mathring{M}$ and $\F U$ is involutive.
\end{proof}

We write $\F_p$ and $\mathcal{E}_p$ for leaves of $\F U$ and $\mathcal{E} U$, and $i_\mathcal{E}$ and $i_\F$ for their inclusions in $M$. Let us then look at these leaves closely.


\begin{lem} \label{souummacaco}
Under the hypothesis of Theorem \ref{ResultNoncompact2}, it holds that $F \circ i_\F (\mathcal{F}_p)$ is an affine space in $\R^{q,n}$ for all $p \in U$. Beyond that, if $q \in i_\mathcal{E}(\mathcal{E}_p)$, then $F \circ i_\F (\F_p)$ and $F \circ i_\F (\F_q)$ are parallel.
\end{lem}

\begin{proof}
Let $p\in U$ be a point and $\mathcal{F}_p$ be the leaf of $\F U$ containing $p$. First, we show this leaf is an affine subspaces of $\R^{q,n}$. Let $q \in i_\F (\mathcal{F}_p) \subset U$ be a point, $f \in \mathcal{F}U_q$ and $X \in T_qM$ be vectors, then equation \eqref{ANormHPnu} implies
\begin{equation} \label{pocoyo2}
A(f,X) = X^j f^i A_{ij} = X^j f^i \frac{1}{\|\H\|}P_{ij}\nu = 0
\end{equation}
because $\F M$ is formed by the vectors in the null space of $P$.

Let $ A_{F \circ i_\F}$ and $A_{i_\F}$ denote the second fundamental tensors of the immersions $F \circ i_{\F}$ and $i_{\F}$.
From equation \eqref{chainrulesecondfund} we know that
$$
A_{F \circ i_\F} = A_F + dF (A_{i_\F}).
$$
On the other hand one can write, for vector fields $f_1,f_2 \in \Gamma(T\F_p)$,
\begin{equation} \label{seconddefinitionAF2}
A_{i_\F}(f_1,f_2) = \nabla_{f_1} f_2 - \nabla'_{f_1} f_2,
\end{equation}
where $\nabla'$ is the Levi-Civita connection of $\F_p$ (with respect to $i_\F$). 
From $\nabla_{f_1} f_2 \in \Gamma(\F U)$ (by equation \eqref{FtimesFinF2}) and $A_{i_\F}(f_1, f_2) \in \Gamma(T\F_p^\perp)$ it holds 
\begin{equation}\label{micuxino}
A_{i_\F} = 0, \hspace{5mm} \mbox{thus} \hspace{5mm} \nabla_{f_1} f_2 =\nabla'_{f_1} f_2,
\end{equation}
Which results in $A_{F \circ i_\F}=0$, i. e. $F \circ i_\F$ is totally geodesic and
\begin{equation}\label{nfdkjflkjhsfdlknfdlksajdlfksj}
D_{f_1} f_2 = A(f_1, f_2) + \nabla_{f_1} f_2 = \nabla_{f_1} f_2 =\nabla'_{f_1} f_2
\end{equation}
implies that the geodesics of $\F_p$ are also geodesics of $\R^{q,n}$, which are straight lines, but $U$ is not geodesically complete, so that the geodesics of $\F_p$ could only be some intervals of these straight lines.

We prove now that $\F_p$ is geodesically complete. Let $\delta: (- a, b) \to \F_p$, $a, b > 0$, be a maximally extended geodesic of $\F_p$ and $\gamma:= i_\F \circ \delta$, as $M$ is geodesically complete, $\gamma$ can be extended $\gamma: \R \to M$, so that $\delta$ could be further extended as long as $\gamma(t) \in U$.


We claim that $\gamma(t) \in U$ for all $t \in [-a,b]$. To prove this we show that $\nabla \|\H\| (\gamma(t)) \neq 0$ for  $t= -a$ and $t=b$. By eqs. \eqref{LaplaceNablaNabla} and \eqref{kernFU} it holds, for every $t \in (-a,b)$,
$$
\nabla_{\dot{\gamma}}\nabla \|\H\| = 0,
$$
so that $\nabla \|\H\|(\gamma(t))$ is the parallel transport of $\nabla \|\H\|(\gamma(0))$. Further, eq. \eqref{pocoyo2} implies
$$
D_{\dot{\gamma}}\nabla \|\H\| = \nabla_{\dot{\gamma}}\nabla \|\H\|
$$
so that $dF(\nabla \|\H\|(\gamma(t)))$ is the parallel transport of $\gamma_0:=dF(\nabla \|\H\|(\gamma(0)))$ over the line\footnote{From eq. \eqref{nfdkjflkjhsfdlknfdlksajdlfksj} $F \circ \gamma$ is a geodesic of $\R^{q,n}$ and thence a straight line.} in $\R^{q,n}$ defined by $dF(\dot{\gamma}(0))$ and $F(\gamma(0))$. This means that $dF(\nabla \|\H\|)(\gamma(t)) = \gamma_0 \neq 0$ for all $t \in [-a,b]$. But $dF$ is linear, so that $\nabla \|\H\|(\gamma(t)) \neq 0$ for all $t \in [-a,b]$, thence $\gamma(t) \in U$ for all $t \in [-a,b]$, which contradicts the maximality of $(-a,b)$. So, $\delta (t)$ is defined for all $t \in \R$ and $F \circ i_\F \circ \delta$ is a whole line in $\R^{q,n}$. Then, analogous to Lemma \ref{jujujujujujjujujujujujjujuj}, $F \circ i_\F(\F_p)$ is an affine $m-r$-dimensional subspace of $\R^{q,n}$.

We fix $p\in U$ and claim that $\F_{p'}$ is parallel to $\F_p$ for any $p' \in \mathcal{E}_p$. As $\mathcal{E}_p$ is a smooth curve, we parametrize it by arc length: $\zeta: [0,a]\to i_{\mathcal{E}}( \mathcal{E}_p)$ with $\zeta(0) =p$ and $\zeta(a)=p'$. Then take $f_p \in \F M_p$ and $f(t)$ the parallel transport of $f_p$ along $\zeta$ with respect to the Levi-Civita connection $\nabla$ of $M$. From $\frac{d}{dt}\langle \dot{\zeta},f(t) \rangle = \langle \nabla_{\dot{\zeta}} \dot{\zeta},f(t) \rangle + \langle \dot{\zeta}, \nabla_{\dot{\zeta}} f(t) \rangle = 0$ follows that $\langle \dot{\zeta},f(t) \rangle = \langle \dot{\zeta}(0),f_p \rangle = 0$ and $f(t) \in \F U$ for all $t \in [0,a]$. 
Then, using eq. \eqref{pocoyo2},
$$
D_{\dot{\zeta}}f(t) = \nabla_{\dot{\zeta}} f(t) + A(\dot{\zeta}, f(t)) = 0,
$$
where $D$ is the Levi-Civita connection of $\R^{q,n}$. The claim follows as in the first case.
\end{proof}

Let us now consider the 1-dimensional leaf of $\mathcal{E}_p$, for some $p \in U$.

\begin{lem}\label{tchaualemanha}
Under the hypothesis of Theorem \ref{ResultNoncompact2}, it holds that the image of $\mathcal{E}_p$ through $F \circ i_{\mathcal{E}_p}$ on $\R^{q,n}$ lies in a plane for every $p \in U$.
\end{lem}

\begin{proof}
$\mathcal{E}_p$ is a smooth curve immersed in $U \subset \mathring{M} \subset M$ through $i_{\mathcal{E}_p}: \mathcal{E}_p \to U$. 
Let $ A_{F \circ i_\mathcal{E}}$ and $A_{i_\mathcal{E}}$ denote the second fundamental tensors of $F \circ i_{\mathcal{E}}$ and $i_{\mathcal{E}}$. From equation \eqref{chainrulesecondfund}:
\begin{equation}\label{patapon}
A_{F \circ i_\mathcal{E}} = A_F + dF (A_{i_\mathcal{E}}).
\end{equation}

On the other hand, for vector fields $e_1, e_2 \in \Gamma(T\mathcal{E}_p)$,
\begin{equation} \label{seconddefinitionAE2}
A_{i_\mathcal{E}}(e_1,e_2) = \nabla_{e_1} e_2 - \nabla'_{e_1} e_2,
\end{equation}

\noindent where $\nabla'$ is the Levi-Civita connection of $\mathcal{E}_p$ (with respect to induced metric). But $\nabla_{e_1} e_2 \in \Gamma(\mathcal{E} M)$ by eq. \eqref{PnablaXYinEU} and $A_{i_\mathcal{E}}(e_1,e_2) \in \Gamma(T\mathcal{E}^\perp)$, so that equation \eqref{seconddefinitionAE2} implies
\begin{equation}\label{ponpata}
A_{i_\mathcal{E}}=0,
\end{equation}
so that $\mathcal{E}_p$ is a geodesic of $U \subset \mathring{M} \subset M$. $\mathcal{E} U_{p'}$ is spanned by $\nabla\|\H\| / \| \nabla\|\H\| \|(p')$ for any $p' \in U$. Let $\gamma (s): (-a,b) \to M$ be the local parametrization by arc lenght of this geodesic in one of the directions $\pm\nabla\|\H\| / \| \nabla\|\H\| \|$ with $\gamma(0)=p$.

We claim that this curve lies in a plane. First of all
\begin{align*}
\frac{d}{ds} (F \circ \gamma) =& dF \circ d\gamma \left( \frac{d}{ds} \right) = \pm dF \left( \frac{\nabla \|\H\|}{\|\nabla \|\H\|\|} \right),\\
\frac{d^2}{ds^2} (F \circ \gamma) 
=& A_{ij} \frac{\nabla^i \|\H\| \nabla^j \|\H\|}{\|\nabla \|\H\|\|^2} = \|\H\| \nu,
\end{align*}

\noindent where we used equation \eqref{ANormHPnu} and $\nabla_{\dot{\gamma}} \dot{\gamma} = 0$. Concerning $\nu$, we get, for any $N \in T_p U^\perp$,
$$
\left< \nabla_{\pm \frac{\nabla \|\H\|}{\|\nabla \|\H\|\|}} \nu, N \right> = 0,
$$
because $\nabla^\perp \nu = 0$. For $f \in \Gamma(\F U)$, we get
\begin{align*}
\left<\nabla_{\pm\frac{\nabla \|\H\|}{\|\nabla \|\H\|\|}} \nu, dF\left(f\right) \right> =& \pm\frac{\nabla \|\H\|}{\|\nabla \|\H\|\|} \left< \nu, dF\left(f\right) \right> - \left< \nu, \pm \nabla_\frac{\nabla \|\H\|}{\|\nabla \|\H\|\|} dF\left(f\right) \right>\\
=&- \frac{1}{\|\H\|}P\left(\frac{\pm\nabla \|\H\|}{\|\nabla \|\H\|\|}, f\right)=0,
\end{align*}
where equation \eqref{ANormHPnu} and the fact that $f$ is in the kernel of $P$ were used. Finally
\begin{align*}
\left< \nabla_{\pm\frac{\nabla \|\H\|}{\|\nabla \|\H\|\|}} \nu, \frac{\pm\nabla \|\H\|}{\|\nabla \|\H\|\|} \right> 
=& -\left< \nu, A\left(\frac{\nabla \|\H\|}{\|\nabla \|\H\|\|}, \frac{\nabla \|\H\|}{\|\nabla \|\H\|\|}\right) + \nabla_{\frac{\nabla \|\H\|}{\|\nabla \|\H\|\|}} \frac{\nabla \|\H\|}{\|\nabla \|\H\|\|} \right>\\
=& -\frac{1}{\|\H\|} P\left(\frac{\nabla \|\H\|}{\|\nabla \|\H\|\|}, \frac{\nabla \|\H\|}{\|\nabla \|\H\|\|}\right) = -\|\H\|,
\end{align*}
so that
\begin{equation}
\nabla_{\pm\frac{\nabla \|\H\|}{\|\nabla \|\H\|\|}} \nu = -\|\H\| dF\left(\pm\frac{\nabla \|\H\|}{\|\nabla \|\H\|\|}\right).
\end{equation}
Let $\mathcal{H}$ be an antiderivative (real or pure imaginary) of $\|\H\|$ restricted to $\gamma$, so that $\dot{\mathcal{H}}(t)=\|\H\|(\gamma(t))$. Then the family of vectorfields $V_\alpha \in \Gamma ((F \circ \gamma)^{-1}(\R^{q,n}) )$, for $\alpha \in \R$ if $\|\H\|^2>0$ or $\alpha= i \beta$ with $\beta \in \R$ if $\|\H\|^2<0$, given by
$$
V_\alpha := \cos (\mathcal{H} + \alpha) dF\left(\pm\frac{\nabla \|\H\|}{\|\nabla \|\H\|\|} \right) - \sin (\mathcal{H} + \alpha) \nu
$$
satisfies $\frac{d}{dt}(V_\alpha) = 0$, this means that any $V_\alpha$ is parallel translated over $F\circ \gamma: \R \to \R^{q,n}$ (thence a constant vector). 
But $\pm dF\left(\frac{\nabla \|\H\|}{\|\nabla \|\H\|\|} \right)$ can be written as a linear combination of two vector fields of the family $V_\alpha$ and lies in the constant plane defined by this family of vector fields and a point of $F\circ \gamma$. 
Thence $F \circ \gamma$ lies in this plane. Plane which is orthogonal to $F \circ i_\F (\F_p)$ because $\frac{\nabla \|\H\|}{\|\nabla \|\H\|\|}$ and $\H$ are orthogonal to any $f \in \F U$.


%
\end{proof}

We want to get a result over the whole manifold $M$, not only on a connected component $U$ of $\mathring{M}$. For that, we need to take a closer look at the set $M \setminus \mathring{M} = \{p \in M : \nabla\|\H\|(p) = 0\}$. In $\overline{\mathring{M}}$ the same equations hold as in $\mathring{M}$, so that we only need to look at the open sets of $M \setminus \mathring{M}$.

As $\mathring{\theta} = \theta$ and $\|\mathring{A}\|^2 = \|A\|^2 - \|P*A\|^2$ (eq. \eqref{AringAPA}) in any open sets of $M \setminus \mathring{M}$, Lemma \ref{mugujininho} is proven in this case exactly as Lemma \ref{voumeembora} in any open subset $V \subset M \setminus \mathring{M}$, thence all equations up to eq. \eqref{srewsrewrewrewerw} of the first case also hold in $V$. Let us consider $V$ maximal, such that $\partial V \subset \partial \mathring{M}$, this implies, for a point $q \in \partial V$, that the tensor $P$ has only one nonzero eigenvalue, and it has to be 1, because on the boundary the equations for $\mathring{M}$ and for $V$ hold (by continuity), but the multiplicity of the eigenvalue 1 in $V$ is constant (because $\nabla P =0$), then $P$ also has only one eigenvector associated with a nonzero eigenvalue in $V$. Beyond this, in $V$, $\|\H\|^2 = \tr(P)=1$ by eq. \eqref{normHsquaredequalr}, so that if there is such a nonempty open set $V$, then $\|\H\|^2>0$ on the whole of $M$, 
because $\|\H\|^2 \neq 0$.

In $V$, one also gets two orthogonal distributions, $\mathcal{E}' V$ and $\mathcal{F}' V$, which are involutive and totally geodesic; beyond this, the leaves of $\mathcal{F}' V$ are affine spaces with $\mathcal{F}'_p \parallel \mathcal{F}'_q$ for any $p,q \in \mathcal{E}'_p \subset V$. In particular, for any $p \in V$, the equations $A_{ij} = P_i^k A_{kj}$ (eq. \eqref{A,PA}) and $\nabla P = 0$ (from Lemma \ref{lemaImportante}) hold. We denote the leaves that contain $p\in V$ by $\mathcal{E}_p'$ and $\mathcal{F}_p'$, and their immersions in $M$ by $i_{\mathcal{E}'}$ and $i_{\mathcal{F}'}$ respectively.

We now prove that Lemma \ref{tchaualemanha} also holds for leaves for the distribution $\mathcal{E}' V$.

\begin{lem}\label{merdadevida}
Under the hypothesis of Theorem \ref{ResultNoncompact2}, it holds that the image of $\mathcal{E}'_p$ through $F \circ i_{\mathcal{E}'_p}$ on $\R^{q,n}$ lies in a plane for every $p \in V$.
\end{lem}

The proof is analogous to the proof of Lemma \ref{tchaualemanha}.

Now let us see what the whole $M$ looks like. We saw then, that the tensor $P$ has globally only one non-zero eigenvector and that the eigenspaces of $P$ give globally the distributions
\begin{align*}
\mathcal{E} M &:= \{V \in M: P(V) = \| \H\|^2 V\}\\
\mathcal{F} M &:= \{V \in M: P(V) = 0\},
\end{align*}
which are involutive and whose leaves ($\mathcal{E}_p$ and $\F_p$) are totally geodesic (by different reasons on $\overline{\mathring{M}}$ and in $M \setminus \overline{\mathring{M}}$). By continuity on the boundary points, all the leaves of $\F M$ are $(m-1)$-dimensional affine spaces of $\R^{q,n}$ and $\F_p ||\F_q$ if $q\in \mathcal{E}_p$.

\begin{lem}\label{fdlksjfdsakljsdfalkjfdaslk}
Under the hypothesis of Theorem \ref{ResultNoncompact2}, let $p \in M$ and $\gamma:\R \to M$ be a, by arc length, parametrization of the leaf $\mathcal{E}_p$, then $\gamma$ lies in a 2-dimensional plane, beyond this the plane is normal to the affine space $F \circ i_\F (\F_q)$, for any $q \in i_\mathcal{E} (\mathcal{E}_p)$.
\end{lem}
\begin{proof}
It follows from Lemas \ref{tchaualemanha} and \ref{merdadevida}.
\end{proof}

Let us see what a particular leaf of $\mathcal{E}_p$ looks like.

\begin{lem}
Under the hypothesis of Theorem \ref{ResultNoncompact2}, there is a $q \in M$ such that $\mathcal{E}_{q}$ is a self-shrinker in $\R^{q,n}$, this means that
$$
\H_{F \circ i_{\mathcal{E}}}(x) = -(F \circ i_{\mathcal{E}}(x))^\perp,
$$
for every $x \in \mathcal{E}_{q_0}$.
\end{lem}
\begin{proof}
Let $q \in M$ be a point such that $\|F(q)\|^2 = \min_{x \in M} {\|F(x)\|^2}$. This implies that
$$
0=\nabla_f \|F\|^2 = \langle \nabla_f F, F \rangle = \langle dF(f), F \rangle,
$$
at $q$, for any $f \in \F M_{q}$.

Let $\delta: \R \to \mathcal{E}_q$ be a, by arc length, parametrization of the leaf $\mathcal{E}_q$ with $i_\mathcal{E}(\delta(0))=q$ and write $\gamma := i_\mathcal{E} \circ \delta$. It holds, for any $q' \in i_\mathcal{E} (\mathcal{E}_q)$, that $F \circ i_\F (\F_{q'}) || F \circ i_\F (\F_{q})$, so that one identifies $f \in \F M_{q} \cong \F M_{q'} \subset \R^{q,n}$ and calculates
\begin{align*}
\left<F \circ \gamma(t) , dF(f) \right> 
&= \int_0^t \left< dF \circ \dot{\gamma}(s), dF(f) \right> ds + \left<F \circ \gamma(0) , dF(f) \right> = 0,
\end{align*}
because $\dot{\gamma}(s) \in \mathcal{E}M \perp \mathcal{F}M$, in particular this means that
\begin{equation}\label{jkjkjkjkjkjkjkjkjkjk}
\langle F \circ i_\mathcal{E} \circ \delta(t), dF(f)\rangle = 0 \,\,\, \forall \, t \in \R.
\end{equation}
Denote $T\mathcal{E}_p^\perp$ the normal bundle of $\mathcal{E}_p$ with respect to $F \circ i_\mathcal{E}$. Then eq. \eqref{jkjkjkjkjkjkjkjkjkjk} implies
$$
\mbox{proj}_{T\mathcal{E}_p^\perp} (F \circ \gamma) = (F \circ \gamma)^\perp.
$$

Otherwise $A_{ij} = \frac{1}{\|\H\|} P_{ij}\nu$ (equation \eqref{ANormHPnu}) in open sets of $\overline{\mathring{M}}$ and $A_{ij} = P_i^k A_{kj}$ in $M \setminus \overline{\mathring{M}}$, so that $A(f,f)=0$ for any $f \in \F M$ and
$$
\H = \mbox{tr}_M A_F = \mbox{tr}_{\mathcal{E}_p} A_F = \mbox{tr}_{\mathcal{E}_p} A_{F\circ i_\mathcal{E}} = \H_{F \circ i_\mathcal{E}},
$$
where we used equations \eqref{popopopopopopopopopopo} and \eqref{AfiequalAfplusDfAi} in the open sets of $M/ \mathring{M}$ and equations \eqref{patapon} and \eqref{ponpata} in $\mathring{M}$ to get $A_F= A_{F \circ i_\mathcal{E}}$.

This implies that
$$
\H_{F \circ i_\mathcal{E}} (\delta(t)) = \H (\gamma(t))= -F^\perp (\delta(t)) = -\mbox{proj}_{T\mathcal{E}_p^\perp} \gamma(t),
$$
so that $i_\mathcal{E}: \mathcal{E}_p \to \R^{q,n}$ is a shrinking self-similar solution of the curve shortening flow.
\end{proof}

Analogously to the first case it holds that $F(M)$ is the product $F(\mathcal{E}_q) \times F(\F_q)$, where $q \in M$ minimizes $\|F\|^2$,
i. e. $F(M)$ is the product of an affine space with a shrinking self-similar solution of the mean curvature flow for plane curves.


\begin{rem}\label{lkjpoijlkj}
At the affine space, the induced (from $\R^{q,n}$) inner product has to be positive definite, because we assumed that $F$ is a spacelike immersion.
\end{rem}
%
%
\begin{rem}
It is not hard to see that $\langle\cdot, \cdot \rangle$ restricted to the plane containing $\mathcal{E}_p$ is positive definite. There one can find a basis made of two orthogonal vectors of length 1, and if one writes the self-shrinking curve in this basis one has just a usual self-shrinker of the curve shortening flow. This is a well studied subject and a classification of such was given by \cite{A.L.86} and can also be found in \cite{Halldorson}. The closed self-shrinkers of the curve shortening flow are called the Abresch \& Langer curves, there are also some curves that ''do not close'' and are dense in some annulus. These curves are not in our classification because they would not satisfy the inverse Lipschitz condition. So that, in our case, the self-shrinking solutions of the mean curvature flow in the plane are just dilatations of the Abresch \& Langer curves in $\E^2$.
\end{rem}
%
\end{proof}
\begin{rem}\label{lolknkjbjklbkjkjnbkjbn}
In particular $M$ is contained in the product of an affine space and a plane, both spacelike, so that $\|\H\|^2>0$.
\end{rem}
\end{section}

\begin{rem}
We found that there are no spacelike self-shrinkers of the mean curvature flow with timelike mean curvature in any of the treated cases, so that Theorems \ref{CompactNonNegative}, \ref{NONCompactNonNegative}, \ref{noMainlyNegative}, \ref{ResultNoncompact1} and Remark \ref{lolknkjbjklbkjkjnbkjbn} sum up to:

\begin{thm}\label{SumUpNoSS}
There are no spacelike self-shrinkers $F:M \to \Rn$ of the MCF that satisfy
\begin{itemize}
 \item $F(M)$ unbounded and $F$ is mainly negative and has bounded geometry or
 \item $\|\H\|^2<0$ and one of the following:
\begin{enumerate}
 \item $M$ is compact.
 \item $F(M)$ is unbounded, $M$ is stochastic complete and\\ $\sup_M \|F\|^2 \leq \infty$.
 \item $F(M)$ is unbounded, $F$ is mainly positive, has bounded geometry and the principal normal parallel in the normal bundle.
\end{enumerate}
\end{itemize}
\end{thm}
\end{rem}

Beyond this, summing up \ref{ResultNoncompact1} and \ref{ResultNoncompact2}, the following classification holds:

\begin{thm}\label{ResultNoncompact3}
Let $M$ be a smooth manifold and $F:M \to \R^{q,n}$ be a mainly positive, spacelike, shrinking self-similar solution of the mean curvature flow with bounded geometry such that $F(M)$ is unbounded. Beyond that, let $F$ satisfy the conditions:  $\|\H\|^2(p) \neq 0$ for all $p \in M$ and the principal normal is parallel in the normal bundle ($\nabla^\perp \nu \equiv 0$). Then one of the two holds:
\begin{align*}
F(M) &= \mathcal{H}_r \times \R^{m-r} \hspace{5mm} \mbox{or}\\
F(M) &= \Gamma \times \R^{m-1},
\end{align*}
where $\mathcal{H}_r$ is an $r$-dimensional minimal surface of the hyperquadric $\mathcal{H}^{n-1}(r)$ (in addition $\|\H\|^2=r>0$) and $\Gamma$ is a rescaling of an Abresch \& Langer curve in a spacelike plane. By $R^{m-r}$ we mean an $m-r$ dimensional spacelike affine space in $\R^{q,n}$.
\end{thm}

\end{section}


\end{document}